\documentclass[11pt]{amsart}
\usepackage{amsmath}
\usepackage{amsfonts}
\usepackage{amssymb}
\usepackage[arrow,matrix,curve,cmtip,ps]{xy}

\numberwithin{equation}{section}

\usepackage{amsthm}

\newtheorem{thm}{Theorem}[section]

\newtheorem{cor}[thm]{Corollary}

\newtheorem{lemma}[thm]{Lemma}
\newtheorem{prop}[thm]{Proposition}
\newtheorem{prob}[thm]{Problem}

\newtheorem{defn}[thm]{Definition}

\theoremstyle{remark}
\newtheorem{remark}[thm]{Remark}

\newcommand{\bb}[1]{\mathbb{#1}}
\newcommand{\cl}[1]{\mathcal{#1}}

\newcommand\id{\mathop{\rm id}}
\newcommand\tr{\mathop{\rm tr}}

\newcommand\diag{\mathop{\rm diag}}

\newcommand\dec{\mathop{\rm dec}}

\newcommand\C{\mathbb{C}}

\newcommand{\OMAX}{\operatorname{OMAX}}
\newcommand{\OMIN}{\operatorname{OMIN}}
\newcommand{\MAX}{\operatorname{MAX}}
\newcommand{\MIN}{\operatorname{MIN}}

\allowdisplaybreaks

\begin{document}

\title[Operator System Structures]{Operator System Structures on Ordered Spaces}

\author[V.~I.~Paulsen]{Vern I.~Paulsen}
\address{Department of Mathematics, University of Houston,
Houston, Texas 77204-3476, U.S.A.}
\email{vern@math.uh.edu}
\author[I.~G.~Todorov]{Ivan G.~Todorov}
\address{Department of Pure Mathematics, Queen's University Belfast,
Belfast BT7 1NN, United Kingdom}
\email{i.todorov@qub.ac.uk}

\author[M.~Tomforde]{Mark Tomforde}
\address{Department of Mathematics, University of Houston, Houston,
  Texas 77204-3476, U.S.A.}
\email{tomforde@math.uh.edu}

\thanks{The first author was supported by NSF Grant DMS-0600191.
The second author was supported by EPSRC Grant D050677/1.
The third author was supported by NSA Grant H98230-09-1-0036.}

\begin{abstract} Given an Archimedean order unit space $(V,V^+,e),$
  we construct a minimal operator system $\OMIN(V)$ and a maximal
  operator system $\OMAX(V)$, which are the analogues of the minimal and
  maximal operator spaces of a normed space. We develop some of the
  key properties of these operator systems and make some progress on
  characterizing when an operator system $\cl S$ is completely
  boundedly isomorphic to either $\OMIN(\cl S)$ or to $\OMAX(\cl S).$
  We then apply these concepts to the study of entanglement breaking
  maps. We prove that for matrix
  algebras a linear map is completely positive from $\OMIN(M_n)$ to
  $\OMAX(M_m)$ if and only if it is entanglement breaking.

\end{abstract}

\maketitle


\section{Introduction}

In the past twenty years, beginning with Ruan's abstract
characterization of operator spaces \cite{Ru}, there has been a great deal of research
activity focused on operator spaces and completely bounded
maps. In contrast, there has been relatively little development of the
abstract theory of operator systems. However, many deep results about operator
spaces are obtained by regarding them as corners of operator systems.
So, potentially, parallel developments in the theory of operator
systems could lead to new insights in the theory of operator spaces.
In this paper we develop the analogues in the operator system setting
of the $\MIN$ and $\MAX$ functors from the category of normed spaces into
the category of operator spaces and study some of their properties.
In particular, we prove that the entanglement breaking maps between
matrix algebras, studied in \cite{Hol}, \cite{HS} and \cite{Ar}, coincide with the
linear maps that are completely positive when the matrix algebra of the
domain is equipped with
our minimal operator system
structure and the target matrix algebra is
equipped with
our maximal operator system structure. This can be
interpreted as showing
that the
entanglement breaking maps are precisely the linear maps between matrix
algebras that are ``universally'' completely positive; i.e., that
remain completely positive independent of the operator system
structures on the domain and the range.

Recall that an operator space is a normed space for which a norm $\|\cdot\|_n$, $n\in \bb{N}$,
is given on the space $M_n(V)$ of $n \times n$ matrices with entries in $V$
in such a way that the family $(\|\cdot\|_n)_{n=1}^{\infty}$ satisfies certain axioms (see \cite[p.181]{Pa} for details).
Given a normed space $X$, there are many inequivalent
operator spaces that all have $X$ as their ``ground level''.
Among all those, there are two distinguished operator spaces,
$\MIN(X)$ and $\MAX(X)$, which represent, respectively, the smallest and
largest operator space structures on $X$. Moreover, $\MIN$ and $\MAX$ can be regarded
as functors from the category whose objects are normed spaces and whose
morphisms are contractive linear maps into the category whose objects are operator spaces
and whose morphisms are completely contractive maps. In this categorical
sense, ``taking the ground level", is really the forgetful functor from the category of
operator spaces to the category of normed spaces, which ignores the structure on the levels above the first.
Effros coined the term \emph{quantization functor} for any functor from the category
of normed spaces into the category of operator spaces for which the forgetful functor
is an inverse. Much work has been done explaining the differences between the $\MIN$ and $\MAX$
functors, constructing other natural operator space structures on normed spaces, and
exploring the behavior of these functors with respect to various natural tensor
norms on each category. These results play a vital role in
the theory of operator spaces and in the theory of $C^*$-algebras.

In this paper we consider a parallel development for operator
systems. Every operator system is at the ground level an ordered
$*$-vector space with an Archimedean order unit and, conversely, given
any Archimedean order unit space, there are possibly many different
operator systems that all have the given Archimedean order unit space
as their ground level. We begin by constructing the analogues of the
$\MIN$ and $\MAX$ functors in this setting, which we denote by $\OMIN$ and $\OMAX$.
Thus, associated with an Archimedean order unit space $V$ we
have two operator systems, $\OMIN(V)$ and $\OMAX(V)$,
whose properties we develop in Section \ref{osys}.
We describe the process of ``Archimedeanization'' of a matrix ordered space
with a matrix order unit, and give a matricial version of the
corresponding result from \cite{pt} concerning ordered $*$-vector spaces.
In Section \ref{mss} we introduce the dual matrix ordered space to a given
matrix ordered space and identify the dual spaces of the operator systems $\OMIN(V)$ and $\OMAX(V)$.
In Section \ref{s_cvs} we provide necessary and sufficient conditions for
an operator system $\cl S$ to be completely boundedly isomorphic to
$\OMIN(\cl S)$ or $\OMAX(\cl S)$.
In Section \ref{s_ent} we apply our results to the study of
entanglement breaking maps encountered in Quantum Information Theory
(see \cite{Ar} and \cite{HS}).
We characterize the entanglement breaking maps
from $M_n$ to $M_k$ as the maps that are completely positive from $\OMIN(M_n)$ to $\OMAX(M_k).$ For maps
between general operator systems we define the notion of a
weak*-entanglement breaking map
and extend some of our results for matrix algebras
to the general setting. The next section is devoted to preliminary notions
and results.

\section{Preliminaries}\label{s_prel}

In this section we recall basic definitions and results and establish
terminology. If $W$ is a real vector space, a {\bf cone} in $W$ is a nonempty
subset $C \subseteq W$ with the following two properties:

\medskip

(a) $\lambda v \in C$ whenever $\lambda \in \bb{R}^+ := [0,\infty)$
and $v \in C$;

(b) $v + w \in C$ whenever $v, w \in C$.

\medskip

\noindent A {\bf $*$-vector space} is a complex vector space $V$
together with a map $^* : V \to V$ which is involutive (i.e., $(v^*)^*
= v$ for all $v \in V$) and conjugate linear (i.e., $(\lambda v + w)^*
= \overline{\lambda} v^* + w^*$ for all $\lambda \in \bb C$ and $v,w \in V$).
If $V$ is a $*$-vector space, then we let $V_h = \{x \in V :
x^* = x \}$ and we call the elements of $V_h$ the {\bf hermitian}
elements of $V$. Note that $V_h$ is a real vector space.

An {\bf ordered $*$-vector space} $(V, V^+)$ is a pair consisting of a $*$-vector space $V$ and a subset $V^+ \subseteq V_h$ satisfying
the following two properties:

\medskip

(a) $V^+$ is a cone in $V_h$;

(b) $V^+ \cap -V^+ = \{ 0 \}$.

\medskip

In any ordered $*$-vector space we may define a partial ordering $\geq$ on $V_h$
by defining $v \geq w$ (or, equivalently, $w \leq v$) if and only if $v - w \in V^+$.
Note that $v \in V^+$ if and only if $v
\geq 0.$ For this reason $V^+$ is called the cone of
{\bf positive} elements of $V$.

If $(V,V^+)$ is an ordered $*$-vector space, an element
$e \in V_h$ is called an {\bf order unit} for $V$ if
for all $v \in V_h$ there exists a real number $r > 0$
such that $re \geq v$.
If $(V, V^+)$ is an ordered $*$-vector space with an order unit $e$,
then we say that $e$ is an {\bf Archimedean order unit} if whenever
$v \in V$ and $re+v \geq 0$ for all real $r >0$, we have that $v \in V^+$.
In this case, we call the triple $(V, V^+, e)$ an {\bf Archimedean
ordered $*$-vector space} or an {\bf AOU space,} for short.

We now recall the notion of a state which will
play a fundamental role in this paper.

\begin{defn} Let $(V, V^+,e)$ be an AOU space and
$s:V \to \bb C$ be a linear map. The map $s$ is called
\emph{unital} if $s(e) =1$, and
\emph{positive} if $s(V^+) \subseteq \bb R^+.$
The map $s$ is called a \emph{state on $V$} if $s$ is unital and positive.
We let $S(V)$ denote
the set of all states on $V$ and call it
the \emph{state space of $V$}.
\end{defn}

\begin{remark} \label{order-norm-rem}
Let $(V, V^+,e)$ be an AOU space.
For $v\in V_h$ let
$$\|v\| = \inf \{ t \in \bb R : -te \le v \le +te \}$$
be the {\em order norm} of $v$ defined in \cite[\S2.2]{pt}.
It was shown in \cite[\S4]{pt} that the order norm $\|\cdot\|$ on $V_h$
can be extended to a
norm on the complex vector space $V,$ but that, in general, this
extension is not unique. Moreover, it was shown in \cite[\S4]{pt} that
among all these extensions there is a
{\em minimal norm} $\|\cdot\|_{m}$ and a
{\em maximal norm} $\|\cdot\|_{M}$, and they satisfy the inequalities
$\|v\|_{m} \le \|v\|_{M} \le 2 \|v\|_{m}$ for every $v \in V$.  (See \cite[Remark~4.8 and Proposition~4.9]{pt}.)
In particular, all extensions of $\|\cdot\|$ to a norm on $V$ are equivalent.
We shall denote by $V_{\min}$ the AOU space $V$ equipped with the minimal norm $\|\cdot\|_{m}$.
It was shown in \cite[\S4.1]{pt} that for $v \in V$ we have
$\|v\|_{m} = \sup \{|s(v)| : s \in S(V)\}.$
\end{remark}

If $V$ is a $*$-vector space, we let $M_{m,n}(V)$ denote the set of
all $m \times n$ matrices with entries in $V$.
The natural addition and scalar multiplication turn $M_{m,n}(V)$
into a complex vector space. We often write $M_{m,n} := M_{m,n}(\bb{C})$, and
let $\{E_{i,j}\}_{i,j=1}^n$ denote its canonical matrix unit system.
If $X = (x_{i,j})_{i,j} \in M_{l,m}$ is a scalar matrix, then for any
$A = (a_{i,j})_{i,j} \in M_{m,n}(V)$ we let $XA$
be the element of $M_{l,n}(V)$ whose $i,j$-entry $(XA)_{i,j}$
equals $\sum_{k=1}^m x_{i,k} a_{k,j}$.
We define multiplication by scalar matrices on the right in a similar way.
Furthermore, when $m=n$, we define a $*$-operation on $M_n(V)$ by letting
$(a_{i,j})_{i,j}^* := (a_{j,i}^*)_{i,j}$. With respect to this operation,
$M_n(V)$ is a $*$-vector space.
We let $M_n(V)_h$ be the set of all hermitian elements of $M_n(V)$.

\begin{defn}
Let $V$ be a $*$-vector space.  We say that $\{ C_n \}_{n=1}^\infty$ is a \emph{matrix ordering} on $V$ if
\begin{enumerate}
\item $C_n$ is a cone in $M_n(V)_h$ for each $n\in \bb{N}$,
\item $C_n \cap -C_n = \{0 \}$ for each $n\in \bb{N}$, and
\item for each $n,m\in \bb{N}$ and each $X \in M_{n,m}(\bb C)$
we have that $X^* C_n X \subseteq C_m$.
\end{enumerate}
In this case we call $(V, \{C_n \}_{n=1}^\infty )$ a \emph{ matrix
  ordered $*$-vector space.} We refer to condition (3) as the
  \emph{compatibility} of the family
  $\{C_n\}_{n=1}^{\infty}$.
\end{defn}

Note that properties (1) and (2) show that $(M_n(V), C_n)$
is an ordered $*$-vector space for each $n\in \bb{N}$. As usual, when $A,B \in M_n(V)_h$, we
write $A \leq B$ if $B-A \in C_n$.

\begin{defn}\label{op-system-abstract-def}
Let $(V, \{ C_n \}_{n=1}^\infty)$ be a matrix ordered $*$-vector space.  For $e \in V_h$ let
$$e_n := \left( \begin{smallmatrix} e & & \\ & \ddots & \\ & & e \end{smallmatrix} \right)$$
be the corresponding diagonal matrix in $M_n(V)$.
We say that $e$ is a \emph{matrix order unit} for $V$
if $e_n$ is an order unit for $(M_n(V), C_n)$ for each $n$.
We say that $e$ is an \emph{Archimedean matrix order unit} if
$e_n$ is an Archimedean order unit for $(M_n(V), C_n)$ for each $n$.
An \emph{(abstract) operator system} is a triple $(V, \{C_n \}_{n=1}^\infty, e)$
where $V$ is a complex $*$-vector space,
$\{C_n \}_{n=1}^\infty$ is matrix ordering on $V$,
and $e \in  V_h$ is an Archimedean matrix order unit.
\end{defn}

We note that the above definition of an operator system was first
introduced by Choi and Effros in \cite{CE2}.  If $V $ and $V'$ are
vector spaces, and $\phi : V \to V'$ is a linear map, then for each
$n\in \bb{N}$ the map $\phi$ induces a linear map $\phi_n : M_n(V)
\to M_n(V')$ given by $\phi_n ((v_{i,j})_{i,j}) :=
(\phi(v_{i,j}))_{i,j}$. If $(V, \{ C_n \}_{n=1}^\infty)$ and $(V',
\{ C_n' \}_{n=1}^\infty)$ are matrix ordered $*$-vector spaces, a
map $\phi : V \to V'$ is called {\bf completely positive} if
$\phi_n(C_n) \subseteq C_n'$ for each $n\in \bb{N}$. Similarly, we
call a linear map $\phi : V \to V'$ a {\bf complete order
isomorphism} if $\phi$ is invertible and both $\phi$ and $\phi^{-1}$
are completely positive.

We denote by $B(\cl H)$ the space of all bounded linear operators
acting on a Hilbert space $\cl H$. A {\bf concrete operator system}
$\cl S$ is a subspace of $B (\cl H)$ such that $\cl S = \cl S^*$ and
$I \in \cl S$. (Here, and in the sequel, we denote by $I$ the
identity operator.) As is the case for many classes of subspaces
(and subalgebras) of $B (\cl H)$, there is an abstract
characterization of concrete operator systems. In this case the
characterization is given by
Definition~\ref{op-system-abstract-def}. If $\cl S \subseteq B (\cl
H)$ is a concrete operator system, then we observe that $\cl S$ is a
$*$-vector space, $\cl S$ inherits an order structure from $B(\cl
H)$, and has $I$ as an Archimedean order unit. Moreover, since $\cl
S \subseteq B(\cl H)$, we have that $M_n(\cl S) \subseteq M_n( B(\cl
H)) \cong B (\cl H^n)$ and hence $M_n(\cl S)$ inherits an order
structure from $B (\cl H^n)$ and the $n \times n$ diagonal matrix
$$\begin{pmatrix} I & & \\ & \ddots & \\ & & I \end{pmatrix}$$ is an
Archimedean order unit for $M_n(\cl S)$.  In other words, $\cl S$ is
an abstract operator system in the sense of Definition
\ref{op-system-abstract-def}. The following result of Choi and
Effros \cite[Theorem~4.4]{CE2} shows that the converse is also true.
For an alternative proof of the result, we refer the reader to
\cite[Theorem~13.1]{Pa}.

\begin{thm}[Choi-Effros]
\label{th_choieffros}
Every concrete operator system $\cl S$ is an operator system.
Conversely, if $(V, \{C_n \}_{n=1}^\infty, e)$ is an operator system, then there
exists a Hilbert space $\cl H$, a concrete operator system
$\cl S \subseteq B(\cl H)$, and a complete order
isomorphism $\phi : V \to \cl S$ with $\phi(e) = I$.
\end{thm}

To avoid excessive notation, we will generally refer to an operator
system as simply a set $V$ with the understanding that $e$ is the
order unit and $M_n(V)^+ := C_n$ is the cone of positive elements in $M_n(V)$.

\section{Operator system structures on AOU spaces}\label{osys}

Let $(V,V^+,e)$ be an AOU space. A {\bf matrix ordering}
on $(V,V^+,e)$ is a matrix ordering $\cl C = \{C_n \}_{n=1}^{\infty}$ on $V$
such that $C_1 = V^+.$
An {\it operator system structure} on $(V, V^+,e)$
is a matrix ordering $\{C_n \}_{n=1}^\infty$ such that
$(V, \{C_n \}_{n=1}^\infty, e)$ is an operator system with $C_1 = V^+.$
Given an operator system $(\cl S, \{ P_n
\}_{n=1}^\infty, e)$ and a unital
positive map $\phi: V \to \cl S$ such that $V^+ = \phi^{-1}(P_1),$ one
obtains an operator system structure on $V$ by setting, $C_n =
\phi_n^{-1}(P_n).$
We shall call this the {\it operator system structure induced by
  $\phi.$} Conversely, given an operator system structure on $V,$ if we let $\cl S = V$ and $\phi$ be the identity map, then we see that the given operator system structure is the one induced by $\phi.$

If $\cl P = \{P_n \}_{n=1}^{\infty}$ and $\cl Q =
\{Q_n \}_{n=1}^{\infty}$ are two matrix orderings on $V$, we say that
$\cl P$ is {\bf stronger} than $\cl Q$ (respectively,
$\cl Q$ is {\bf weaker} than $\cl P$) if $P_n\subseteq Q_n$ for all
$n\in \bb{N}$. Note that $\cl P$ is stronger than $\cl Q$ if and only
if for every $n,$ and every $A,B \in M_n(V)_h,$ the inequality
$A \le_{\cl P} B$ implies that $A \le_{\cl Q} B,$
where the subscripts are used to denote
the partial orders induced by $\cl P$ and $\cl Q,$ respectively. Equivalently,
$\cl P$ is stronger than $\cl Q$ if and only if the identity map on $V$ is
completely positive from $(V, \{P_n \}_{n=1}^\infty)$ to $(V, \{Q_n \}_{n=1}^\infty).$

In this section we wish to describe the various operator system
structures with which an AOU space can be equipped.
We shall prove that every AOU space possesses a strongest and a weakest
operator system structure, which we will call the {\bf maximal} and {\bf minimal}
operator system structure, respectively, and we shall
characterize the corresponding matrix orderings.

We begin with the weakest operator system structure.
By  Kadison's Representation Theorem (see \cite{ka} and \cite[Theorem~II.1.8]{Alf}, and also see
\cite[Theorem~5.2]{pt} for the precise
statement which we shall use),
given an AOU space $(V, V^+,e)$ there exists a
compact Hausdorff topology on $S(V)$
such that the unital linear map $\Phi: V \to C(S(V))$
into the $C^*$-algebra of continuous functions on $S(V)$ defined by
$\Phi(v)(s)= s(v)$ is an order isomorphism onto its range.
Equivalently, we have $V^+ = \Phi^{-1}(P_1),$ where $P_1$ denotes the set of
non-negative valued continuous functions on $S(V).$
Since unital $C^*$-algebras are operator systems, Kadison's Representation Theorem induces an operator system
structure $\{ C_n \}_{n=1}^\infty $ on $V.$ We have that $(v_{i,j})
\in C_n$ if and only if $(\Phi(v_{i,j})) \ge 0$ in $M_n(C(S(V)),$
if and only if $(s(v_{i,j})) \in M_n^+,$ for every $s \in S(V).$
We shall prove that this operator system structure is the desired weakest
operator system structure, that is, the one for which the cones of
positive elements are as large as possible.

\begin{defn} \label{C-n-min-def}
Let $(V, V^+,e)$ be an AOU space.  For each $n \in \bb N$
set
$$C_n^{\min}(V) = \left\{ (v_{i,j}) \in M_n(V) :
  \sum_{i,j=1}^n \overline{\lambda}_i \lambda_j v_{i,j} \in V^+ \text{
    for all } \lambda_1, \ldots, \lambda_n \in \bb C \right\}$$
and let $\cl C^{\min}(V) = \{ C_n^{\min} \}_{n=1}^\infty$.
\end{defn}

The following result gives an alternative way to describe the elements
of $\cl C^{\min}(V)$, while simultaneously proving that
it is an operator system structure on $V$.

\begin{thm} \label{CnMin-states}
Let $(V,V^+,e)$ be an AOU space and $n\in \bb{N}$.
Then $(v_{i,j}) \in C_n^{\min}(V)$ if and only if
$(s(v_{i,j})) \in M_n^+$ for each $s \in S(V).$ Hence $\cl
C^{\min}(V)$ is the operator system structure on $V$ induced by the inclusion of
$V$ into $C(S(V)).$
\end{thm}
\begin{proof}
Suppose $(v_{i,j}) \in M_n(V)$.  Then we see that
\begin{align*}
&(s(v_{i,j})) \in M_n^+ \text{ for all } s \in S(V) \\
\iff &\langle (s(v_{i,j})) \vec{x}, \vec{x} \rangle \geq 0 \text{
  for all $s: V \to \bb C$ and all } \vec{x} =
\left( \begin{smallmatrix} \lambda_1 \\ \vdots \\
    \lambda_n \end{smallmatrix} \right) \in \bb C^n \\
\iff  &\sum_{i,j=1}^n \overline{\lambda}_i \lambda_j s(v_{i,j}) \geq 0
\text{ for all $s: V \to \bb C$ and all }
\left( \begin{smallmatrix} \lambda_1 \\ \vdots \\
    \lambda_n \end{smallmatrix} \right) \in \bb C^n \\
\iff &s \left( \sum_{i,j=1}^n \overline{\lambda}_i \lambda_j v_{i,j}
\right) \geq 0  \text{ for all $s: V \to \bb C$ and all }
\left( \begin{smallmatrix} \lambda_1 \\ \vdots \\
        \lambda_n \end{smallmatrix} \right) \in \bb C^n \\
\iff &\sum_{i,j=1}^n \overline{\lambda}_i \lambda_j v_{i,j} \in V^+
\text{ for all } \left( \begin{smallmatrix} \lambda_1 \\ \vdots \\
    \lambda_n \end{smallmatrix} \right) \in \bb C^n \quad \text{ (by \cite[Proposition 3.13]{pt})}\\
\iff &(v_{i,j}) \in C_n^{\min}(V).
\end{align*}
\end{proof}

\begin{defn} \label{min-operator-system-def}
Let $(V,V^+,e)$ be an AOU space.  We define $\OMIN(V)$ to be the
operator system $(V,\cl C^{\min}(V),e)$.
\end{defn}

Thus, up to complete order isomorphism, $\OMIN(V)$ can be identified
with a subspace of $C(S(V)).$
We next examine its universal properties.

\begin{thm}\label{autocp}
Let $(V, V^+,e)$ be an AOU space.
If $(W, \{C_n\}_{n=1}^\infty)$ is a matrix ordered
$*$-vector space and $\phi : W \to \OMIN(V)$ is a positive
linear map, then $\phi$ is completely positive.

Moreover, if $\tilde{V} = (V,\{\tilde{C}_n\}_{n=1}^\infty,e)$ is an operator
system with $\tilde{C}_1 = V^+$ and such that for every operator system $W$,
any positive map $\psi : W\rightarrow\tilde{V}$ is completely
positive, then the identity map is a unital complete order isomorphism
between $\tilde{V}$ and $\OMIN(V)$.
\end{thm}
\begin{proof}
(i) Let $A = (a_{i,j}) \in C_n$.  If
$X = \left( \begin{smallmatrix} \lambda_1 \\ \vdots
    \\ \lambda_n \end{smallmatrix} \right) \in M_{n,1}(\bb C)$ then
$X^* AX = \sum_{i,j=1}^n \overline{\lambda}_i \lambda_j a_{i,j} \in
C_1$.  Since $\phi$ is positive, this implies that
$$\sum_{i,j=1}^n \overline{\lambda}_i \lambda_j \phi(a_{i,j}) =
\phi \left(  \sum_{i,j=1}^n \overline{\lambda}_i \lambda_j a_{i,j} \right) \in V^+.$$
Thus $(\phi(a_{i,j})) \in C_n^{\min}(V)$, and $\phi$ is completely positive.

(ii) Let $\iota : \tilde{V}\rightarrow \OMIN(V)$ be the identity map. By
(i), $\iota$ is completely positive, and by the assumption,
$\iota^{-1}$ is completely positive. Since $\iota$ is also unital,
we have that $\tilde{V}$ and $\OMIN(V)$ are completely order
isomorphic.
\end{proof}

\begin{cor}\label{c_weakestos}
Let $(V, V^+, e)$ be an AOU space.
If $(V, \{C_n \}_{n=1}^\infty, e)$
is any operator system structure on $V$ then $C_n \subseteq C_n^{\min}(V)$
for all $n\in \bb{N}$.
\end{cor}
\begin{proof}  The identity map from  $(V, \{C_n \}_{n=1}^\infty, e)$
  to $\OMIN(V)$ is positive and, hence, completely positive by
  Theorem \ref{autocp}. Thus, $C_n\subseteq C_n^{\min}(V)$ for each $n\in \bb{N}$.
\end{proof}

Corollary \ref{c_weakestos} shows that
$\OMIN(V)$ is the weakest operator system structure that an AOU space
$V$ can be equipped with.

We note that, by virtue of Theorem \ref{th_choieffros},
every operator system is also an operator space.
We next identify the operator space structure of $\OMIN(V).$

\begin{prop}\label{p_minomin}
Let $(V,V^+,e)$ be an AOU space and $V_{\min}$ denote the vector space
$V$ equipped with the minimal norm $\| \cdot \|_m$
(see Remark~\ref{order-norm-rem}).  Then the identity map on $V$ is a complete
isometry between the operator spaces $\MIN(V_{\min})$ and $\OMIN(V).$
\end{prop}
\begin{proof}
The canonical inclusion $V \to \widehat{V} \subset C(S(V))$
is an isometry on $V_{\min}$ and a complete order isomorphism between
$\OMIN(V)$ and the subspace $\widehat{V}$ together with the
operator system structure that it inherits.
Hence, for $(v_{i,j}) \in M_n(V),$ we have that
$$\|(v_{i,j})\|_{M_n(\OMIN(V))} = \|(\hat{v}_{i,j})\|_{M_n(C(S(V)))}
 \ge \|(v_{i,j})\|_{M_n(V_{\min})},$$
with the last inequality following from the fact that the map $v \to \hat{v}$
is an isometry on $V_{\min}$ and the fact that $\MIN(V_{\min})$
is the smallest of all possible operator space structures.

However, each state $s \in S(V)$ is a contractive
linear functional on $V_{\min},$ and hence completely
contractive on $\MIN(V_{\min})$; thus,
$$\|(v_{i,j})\|_{M_n(\OMIN(V))} = \sup\mbox{}_{s \in S(V)}\|(s(v_{i,j}))\|_{M_n}
\le \|(v_{i,j})\|_{M_n(\MIN(V_{\min}))},$$
and the result follows.
\end{proof}

We now turn our attention to the maximal operator system structure on
an AOU space.
Given a $*$-vector space $V,$ we identify the vector space $M_n(V)$ of all $n \times n$ matrices
with entries in $V$ with the (algebraic) tensor product
$M_n\otimes V$ in the natural way.

\begin{lemma}\label{hermat}
We have that $M_n(V)_h = (M_n)_h \otimes V_h$ (the right hand side
being the algebraic tensor product of real vector spaces).
\end{lemma}
\begin{proof} It is obvious that $(M_n)_h \otimes V_h\subseteq M_n(V)_h$.
Conversely, suppose that $v = (v_{i,j})\in M_n(V)_h$. Then
$v_{i,j}^* = v_{j,i}$, $i,j = 1,\dots,n$. Write $v = \sum_{i=1}^n
E_{i,i}\otimes v_{i,i} + \sum_{i < j} (E_{i,j}\otimes v_{i,j} +
E_{j,i}\otimes v_{j,i})$. Clearly, $\sum_{i=1}^n E_{i,i} \otimes
v_{i,i} \in (M_n)_h \otimes V_h$. Fix $i,j$ with $i < j$ and write
$v_{i,j} = \sum_{k=1}^4 \lambda_k w_k$, where $w_k\in V^+$ and
$\lambda_k\in \bb{C}$. Then $v_{j,i} = v_{i,j}^* = \sum_{k=1}^4
\overline{\lambda_k} w_k$ and so
\begin{eqnarray*}
E_{i,j}\otimes v_{i,j} + E_{j,i}\otimes v_{j,i} & = & E_{i,j}\otimes
\left(\sum_{k=1}^4 \lambda_k w_k\right) + E_{j,i}\otimes
\left(\sum_{k=1}^4 \overline{\lambda_k} w_k\right)\\ & = &
\sum_{k=1}^4 \left((\lambda_k E_{i,j})\otimes w_k +
(\overline{\lambda_k}E_{j,i})\otimes w_k\right)\\ & = &
\sum_{k=1}^4 (\lambda_k E_{i,j} +
\overline{\lambda_k}E_{j,i})\otimes w_k\in (M_n)_h \otimes V_h.
\end{eqnarray*}
It follows that $v\in (M_n)_h \otimes V_h$.
\end{proof}

\begin{defn}\label{d_min}
Let $(V,V^+)$ be an ordered $*$-vector space.
Set
$$D_n^{\max}(V) = \left\{\sum_{i=1}^k a_i\otimes v_i : v_i\in V^+, a_i\in
M_n^+, i = 1,\dots,k, \ k\in\bb{N}\right\}$$ and
$\cl D^{\max}(V) = \{D_n^{\max}(V)\}_{n=1}^{\infty}.$
\end{defn}

\begin{lemma}\label{min}
Let $(V,V^+)$ be an ordered $*$-vector space. Suppose that
$P_n\subseteq M_n(V)_h$ is a cone for each $n\in\bb{N}$, the family
$\{P_n\}_{n=1}^{\infty}$ is a compatible matrix ordering, and $P_1 =
V^+$. Then $D_n^{\max}(V) \subseteq P_n$, for each $n\in \bb{N}$.
\end{lemma}
\proof
Suppose that $\{P_n\}_{n=1}^{\infty}$ is a compatible collection of
cones with $P_1 = V^+$. If $X\in M_{n,1}$ then
$$X V^+ X^* = X P_1 X^* \subseteq P_n.$$ It
follows that $a\otimes v\in P_n$ for each $a\in M_n^+$ of rank
one and each $v\in V^+$.
Since every element of $M_n^+$ is the sum of rank one elements of
$M_n^+$, we conclude that $a\otimes v\in P_n$ for
all $a\in M_n^+$ and all $v\in V^+$. Thus $D_n^{\max}(V) \subseteq
P_n$ for all $n \in \mathbb{N}$.
\endproof

If $v_1,\dots,v_m\in V$ we let
$\diag(v_1,\dots,v_m)$
denote the element of $M_m(V)$ with $v_1,\dots,v_m$ on its diagonal
(in this order) and zeros elsewhere.

\begin{prop}\label{minmax}
Let $(V,V^+,e)$ be an AOU space. The
following hold:
\begin{itemize}
\item[(i)] $\cl D^{\max}(V)$ is a matrix ordering on $V$ and $e$ is
a matrix order unit for this ordering;

\item[(ii)]  $D_n^{\max}(V) = \{\alpha \diag(v_1,\dots,v_m)\alpha^* :
\alpha\in M_{n,m}, v_i\in V^+, i= 1,\dots,m,$ $m\in\bb{N}\}$;

\item[(iii)]  if $\{P_n\}_{n=1}^{\infty}$ is a matrix ordering on $V$ with $P_1 =
V^+$ then
$D_n^{\max}(V)\subseteq P_n\subseteq C_n^{\min}(V)$ and $e$ is a matrix
order unit for $\{P_n\}_{n=1}^\infty.$
\end{itemize}
\end{prop}
\begin{proof} Let $D_n$ denote the right hand side of the equation in (ii). We first
observe that $D_n$ is a cone in $M_n(V)_h$. If $v_1,\dots,v_m\in
V^+$ and $\alpha = (\alpha_{i k})\in M_{n,m}$ then the
$(i,j)$-entry of $\alpha\diag(v_1,\dots,v_m)\alpha^*$ is equal to
$\sum_{k=1}^m \alpha_{i k}\overline{\alpha_{j k}} v_k$, and hence
$$\left(\sum_{k=1}^m \alpha_{i k}\overline{\alpha_{j k}} v_k\right)^* =
\sum_{k=1}^m \alpha_{j k}\overline{\alpha_{i k}} v_k.$$ This shows
that $D_n\subseteq M_n(V)_h$. It is obvious that $D_n$ is closed
under taking multiples with non-negative real numbers.
If
$$\alpha \diag(v_1,\dots,v_m)\alpha^*, \
\beta \diag(w_1,\dots,w_k)\beta^*\in D_n$$
then
\begin{eqnarray*}
& & \alpha \diag(v_1,\dots,v_m)\alpha^* + \beta
\diag(w_1,\dots,w_k)\beta^*\\ & = & [\alpha \ \ \beta]
\diag(v_1,\dots,v_m,w_1,\dots,w_k) [\alpha \ \ \beta]^*\in D_n \ ;
\end{eqnarray*}
in other words, $D_n$ is a cone. If $\alpha
\diag(v_1,\dots,v_m)\alpha^*\in D_n$ and $\beta\in M_{k,n}$ then
$$\beta(\alpha \diag(v_1,\dots,v_m)\alpha^*)\beta^* =
(\beta\alpha)\diag(v_1,\dots,v_m)(\beta\alpha)^*\in D_k,$$ so
$\{D_n\}_{n=1}^{\infty}$ is compatible. It is clear that $D_1 =
V^+$. Lemma \ref{min} now implies that $D_n^{\max}\subseteq D_n$
for each $n\in \bb{N}$.

On the other hand, if $v_1,\dots,v_m\in V^+$ then
$\diag(v_1,\dots,v_m) = \sum_{i=1}^m E_{i,i}\otimes v_i\in
D_m^{\max}$. By the compatibility of $(D_n^{\max})_{n=1}^{\infty}$,
for any $\alpha \in M_{n,m}$ we have that
$$\alpha\diag(v_1,\dots,v_m)\alpha^*\in D_n^{\max}$$
for every $\alpha\in M_{n,m}$.
Thus, $D_n\subseteq D_n^{\max}$ and (ii) is established.

We claim that $e_n$ is an order unit for $D_n^{\max}$. Let $w
\in M_n(V)_h$. Assume first that $w = a\otimes v$, where $a\in
(M_n)_h$ and $v\in V_h$. Write $a = a_1 - a_2$ and $v = v_1 -
v_2$, where $a_i\in M_n^+$ and $v_i \in V^+$, $i = 1,2$. Then
$$a\otimes v = a_1\otimes v_1 - a_1\otimes v_2 - a_2\otimes v_1 +
a_2\otimes v_2.$$ Let $r,s\in \bb{R}^+$ be such that
$$-rI_n \leq a_i \leq rI_n, \ \ se \pm v_i\in V^+, \quad i = 1,2.$$ Then
$rs e_n - (\pm a_i\otimes v_j)\in D_n^{\max}$ and so $4rs e_n - w
\in D_n^{\max}$.

Now let $w\in M_n(V)_h$ be arbitrary. By Lemma \ref{hermat}, $w =
\sum_{j=1}^N a_j\otimes v_j$, where $a_j\in (M_n)_h$ and $v_j\in
V_h$, $j = 1,\dots,N$. The claim now follows from the previous
paragraph.

Suppose that $\cl P =\{P_n\}_{n=1}^{\infty}$ is a
matrix ordering on $(V,V^+,e)$.
If we let $W = (V,\cl P)$ then, by
Theorem~\ref{autocp}, the identity map from $W$ to $\OMIN(V)$ is
completely positive and hence $\cl P$ is stronger than $\cl
C^{\min}(V)$; in other words, $P_n \subseteq C_n^{\min},$ for all $n\in\bb{N}$.
The inclusion $D_n^{\max} \subseteq P_n$ follows from Lemma~\ref{min}.
In particular, $D_n^{\max} \subseteq C_n^{\min}$, and hence
$D_n^{\max} \cap -D_n^{\max}\subseteq C_n^{\min}\cap -C_n^{\min} = \{0\}$,
$n\in \bb{N}$,
which shows that $\cl D^{\max}(V)$ is a matrix ordering.

Finally, if $A \in M_n(V)_h,$ then there exists $r>0,$ so that $re_n
+A \in D_n^{\max}$ and hence, $re_n +A \in P_n,$ so that $e_n$ is an
matrix order unit for $\cl P.$
\end{proof}

\begin{remark}\label{r_two}
Inspection of the above proof shows that Proposition~\ref{minmax}(ii),
as well as the inclusion $D_n^{\max}(V)\subseteq P_n$,
hold for matrix ordered $*$-vector spaces
not necessarily possessing an order unit.  Furthermore, the above results show that
 among all matrix orderings on
  an AOU space $(V, V^+,e)$,  the matrix ordering $\cl D^{\max}(V)$ is the strongest while
  $\cl C^{\min}(V)$ is the weakest. Also, $\cl C^{\min}(V)$ is
  simultaneously the
  weakest among all operator system structures on $(V, V^+, e).$
If $e$ was an Archimedean matrix order unit for the matrix ordering
$\cl D^{\max}(V),$ then we would also have that $\cl D^{\max}(V)$ is the
strongest operator system structure on $V.$ Unfortunately, this is not
generally the case and we sketch in the details of an example in the
next remark. Thus, generally, $(V, \cl D^{\max}(V), e)$ is not an
operator system.
For this reason we need to discuss the Archimedeanization process for matrix
ordered spaces. This theory was developed in detail for ordered $*$-vector
spaces in \cite[\S 2.3 and \S 3.2]{pt}.
\end{remark}

\begin{remark} \label{Dnotarch}  Let $V= C([0,1])$ denote the vector space of
  continuous complex-valued functions on the unit interval, with $V^+$
  the usual cone of positive functions and $e$ the constant
  function taking value 1, and
let $P(t) =  \begin{pmatrix} 1 & e^{2 \pi it} \\ e^{-2 \pi it} &
  1 \end{pmatrix} \in M_2(V)_h.$
It can be shown by a rather long calculation that for every $r>0,$ we
have that $re_2 + P(t) = \begin{pmatrix} 1+r & e^{2\pi it}\\ e^{-2 \pi it} &
  1+r \end{pmatrix} \in D^{\max}_2(V)$, and
  $P(t) \notin D^{\max}_2(V).$ This shows that $e=1$ is not an
  Archimedean matrix order unit for the matrix ordering $\cl D^{\max}(V).$
  We sketch
the proof of these claims.

To see the second claim, one assumes that $P(t) = \sum_{j=1}^m Q_j
\otimes p_j(t)$ with $Q_j= \begin{pmatrix} a_j & b_j \\ \overline{b_j} & c_j \end{pmatrix}\in M_2^+$ and $p_j(t) \in V^+,$ for all
$j.$ Thus, $1 = \sum_{j=1}^m a_jp_j(t) = \sum_{j=1}^m c_jp_j(t),$
while $e^{2 \pi it} = \sum_{j=1}^m b_j p_j(t).$ One uses the fact that
$|b_j|^2 \le a_jc_j,$ to show that this set of equalities is
impossible.

To show that $re_2 + P \in D^{\max}_2(V),$ for every $r>0,$ one first
shows that for every $\epsilon > 0,$ there exists $Q \in
D^{\max}_2(V),$ with $||P(t) - Q(t)||_{M_2} < \epsilon$ for all $t.$
Then one shows that $H \in M_2(V)_h$ and $||H(t)|| < \epsilon$ for all $t,$
implies that $4 \epsilon e_2 + H \in D^{\max}_2(V).$ Hence, $re_2 + P =
[re_2 + (P-Q)] + Q \in D^{\max}_2(V),$ where $Q$ is chosen as above for
$\epsilon = r/4.$
\end{remark}

We describe the Archimedeanization process for matrix ordered spaces
below in somewhat more detail than is needed for the special case
of $\cl D^{\max}(V)$ since the general results are likely to
be useful for other situations as well.

\subsection{The Archimedeanization of a matrix ordered $*$-vector space with a matrix order unit} \label{Arch-Operator-Systems-sec}

It was shown in \cite[\S 3.2]{pt} that for any ordered $*$-vector space $(V,V^+)$ with order unit $e$,
there is a functorial way to produce an AOU space, called the \emph{Archimedeanization of $V$},
which is the largest quotient of $V$ containing the class of $e$ as an Archimedean order unit.
Specifically, if $(V,V^+,e)$ is an ordered $*$-vector space with order unit $e$, we define
$D:= \{ v \in V_h : re + v \in V^+ \text{ for all $r > 0$} \}$ and $N := \displaystyle\bigcap_{f\in S(V)} \ker f$.
Then $N$ is a complex subspace of $V$ closed under the $*$-operation, so that the quotient $V/N$ is a $*$-vector space
in the natural way and $(V/N)_h = \{ v+N : v \in V_h \}$.  We define an order on $V/N$ by $(V/N)^+ := \{ v +N : v \in D \}$.
The Archimedeanization of $V$ is defined to be $V_\textnormal{Arch} := (V/N, (V/N)^+, e)$, and it is shown in
\cite[Theorem~3.16]{pt} that $V_\textnormal{Arch}$ is an AOU space that is characterized by the following universal property:
the quotient $q : V \to V_\textnormal{Arch}$ is a positive linear map and whenever $(W,W^+, e')$ is an AOU space and
$\phi : V \to W$ is a unital positive linear map, there exists a unique positive linear map
$\widetilde{\phi} : V_\textnormal{Arch} \to W$ with $\phi = \widetilde{\phi} \circ q$.
\begin{equation*}
\xymatrix{  V_\textnormal{Arch} \ar@{-->}[rd]^{\widetilde{\phi}} & \\ V \ar[u]^q  \ar[r]^<>(.35)\phi & W}
\end{equation*}

Here we shall generalize this construction to matrix ordered $*$-vector spaces containing a matrix order unit.

\begin{defn}
Let $(V, \{C_n \}_{n=1}^\infty)$ be a matrix ordered $*$-vector space with matrix order unit $e$.  For each $n \in \mathbb{N}$ define $N_n := \displaystyle\bigcap_{f \in S(M_n(V))} \ker f$.  Note that using the notation of the previous paragraph, we have $N = N_1$.
\end{defn}

\begin{lemma} \label{N-matrix-is-Nn-lem}
If $(V, \{C_n \}_{n=1}^\infty)$ is a matrix ordered $*$-vector space with matrix order unit $e$, then for each $n \in \mathbb{N}$ we have $N_n = M_n (N)$.
\end{lemma}

\begin{proof}
If $A = (a_{k,l})_{k,l} \in N_n$, then every state on $M_n(V)$ annihilates $A$,
and consequently every positive linear functional on $M_n(V)$ annihilates $A$.
If $s : V \to \C$ is any state on $V$ and $P = (p_{k,l})_{k,l} \in M_n^+$ is any positive matrix over $\C$, then the map $\tilde{s}_P : M_n(V) \to \C$ given by $\tilde{s}_P ((x_{k,l})_{k,l}) := \sum_{k,l} s(p_{k,l} x_{k,l})$ is a linear functional on $M_n(V)$.  Furthermore, we can argue that $\tilde{s}_P$ is positive as follows: Any rank one positive matrix $P \in M_n$ has the form $P = \alpha^* \alpha$ for $\alpha \in M_{1,n}$, and for any $X=(x_{k,l})_{k,l} \in C_n$, we have that $\alpha^* X \alpha \in C_1$.  Hence $\tilde{s}_P ((x_{k,l})_{k,l}) = \sum_{k,l} s(p_{k,l} x_{k,l}) = s( \sum_{k,l} \alpha_k x_{k,l} \overline{\alpha_l} ) = s (\alpha^* X \alpha) \geq 0$.  Since any positive matrix $P \in M_n$ is the sum of rank one positive matrices, the linearity of $s$ shows that $\tilde{s}_P ((x_{k,l})_{k,l}) \geq 0$ for all $P \in M_n^+$ and all $(x_{k,l})_{k,l} \in C_n$.

It follows that $\tilde{s} (A) = 0$, and $s ( \sum_{k,l} p_{k,l} a_{k,l} ) = 0$
for every state $s : V \to \C$ and every positive matrix $P = (p_{k,l})_{k,l} \in M_n^+$.
If we choose $1 \leq k \leq n$, and let $D$ be the diagonal matrix
with a $1$ in the $(k,k)$ position and zeroes elsewhere,
then $D \in M_n^+$ and $\tilde{s}_D(A) = s(a_{k,k}) = 0$ so that
\begin{equation}\label{diag-zero-eqn}
s(a_{k,k}) = 0 \qquad \text{for all states $s : V \to \C$.}
\end{equation}
Furthermore, if we choose $1 \leq k, l \leq n$ and let $\alpha \in M_{1,n}$ be the row vector with a $1$ in the $k$\textsuperscript{th} position, a $1$ in the $l$\textsuperscript{th} position, and zeroes elsewhere, then $P := \alpha^* \alpha \in M_n^+$.  Since $P$ has 1's in $(k,k)$, $(k,l)$, $(l,k)$, and $(l,l)$ positions, and zeroes elsewhere, we see that $\tilde{s}_P (A) = s (a_{k,k}) + s (a_{k,l}) + s (a_{l,k}) + s (a_{l,l}) = 0$.  Using \eqref{diag-zero-eqn} we see that
\begin{equation} \label{off-diags-zero-eqn}
s (a_{k,l}) + s (a_{l,k}) = 0  \qquad \text{for all states $s : V \to \C$.}
\end{equation}
Similarly, if we let $\beta \in M_{1,n}$ be the row vector with a $1$ in the $k$\textsuperscript{th} position, an $i$ in the $l$\textsuperscript{th} position, and zeroes elsewhere, then $Q := \beta^* \beta \in M_n^+$.
Since $Q$ has 1 in $(k,k)$ and $(l,l)$ positions, $i$ in the $(k,l)$ position, $-i$ in the $(l,k)$ position, and zeroes elsewhere, we see that $\tilde{s}_P (A) = s (a_{k,k}) + is (a_{k,l}) -i s (a_{l,k}) + s (a_{l,l}) = 0$.  Using \eqref{off-diags-zero-eqn} we see that
\begin{equation} \label{off-diags-comp-zero-eqn}
is (a_{k,l}) -is (a_{l,k}) = 0  \qquad \text{for all states $s : V \to \C$.}
\end{equation}
It follows from \eqref{off-diags-zero-eqn} and \eqref{off-diags-comp-zero-eqn} that for any $1 \leq k, l \leq n$ it is the case that $s (a_{k,l}) = 0$ for all states $s : V \to \C$.  Thus $a_{k,l} \in N$ and $A \in M_n(N)$.  Hence $N_n \subseteq M_n(N)$.

For the converse, suppose
that $A = (a_{k,l})_{k,l}  \in M_n(N)$.  Let $s : M_n(V) \to \C$ be a state on $M_n(V)$.
For $1 \leq k,l \leq n$ define $s_{k,l} : V \to \C$ by $s_{k,l} (v) := s (E_{k,l} \otimes v)$,
where $E_{k,l}$ is the matrix with a $1$ in the $(k,l)$ position and zeroes elsewhere.
Then the $s_{k,l}$'s are linear functionals and $s(A) = \sum_{k,l} s_{k,l}(a_{k,l})$.
Choose $1 \leq k \leq n$.  For any $v \in V^+ = C_1$ we see that the matrix
$\diag (0, \ldots, 0, v, 0, \ldots 0)$, with a $v$ in the $(k,k)$ position is
positive since $\diag (0, \ldots, 0, v, 0, \ldots 0) =
(0, \ldots, 0, 1, 0, \ldots 0)^* (v) (0, \ldots, 0, 1, 0, \ldots 0) \in C_n$.
Thus $s_{k,k}(v) = s( \diag (0, \ldots, 0, v, 0, \ldots 0)) \geq 0$.  Hence $s_{k,k} : V \to \C$ is a
positive linear functional.  It follows that
\begin{equation} \label{diag-plf}
s_{k,k}(x) = 0 \qquad \text{ for all $x \in N$.}
\end{equation}
In addition, if $1 \leq k,l \leq n$ and $v \in V^+ = C_1$, then the matrix $P \in M_n(V)$ with $v$ in the $(k,k)$, $(k,l)$, $(l,k)$, and $(l,l)$ positions, and zeroes elsewhere, is an element of $C_n$.  It follows that $s_{k,k} (v) + s_{k,l}(v) + s_{l,k}(v) + s_{l,l}(v) = s(P) \geq 0$.  Thus $s_{k,k} + s_{k,l} + s_{l,k} + s_{l,l}$ is a positive linear functional.  It follows that $s_{k,k}(x) + s_{k,l}(x) + s_{l,k}(x) + s_{l,l}(x) = 0$ for all $x \in N$, and using \eqref{diag-plf} we have that
\begin{equation} \label{off-diag-plf-1}
s_{k,l}(x) + s_{l,k}(x) = 0 \text{ for all $x \in N$.}
\end{equation}
Similarly, fix $1 \leq k,l \leq n$.  If $v \in V^+ = C_1$, then the matrix $Q \in M_n(V)$ with $v$ in the $(k,k)$ and $(l,l)$ position, $iv$ in the $(k,l)$ position, $-iv$ in the $(l,k)$ position, and zeroes elsewhere, is an element of $C_n$.  It follows that $s_{k,k} (v) + is_{k,l}(v) -i s_{l,k}(v) + s_{l,l}(v) = s(P) \geq 0$.  Thus $s_{k,k} + is_{k,l} - is_{l,k} + s_{l,l}$ is a positive linear functional.  It follows that $s_{k,k}(x) + is_{k,l}(x) - is_{l,k}(x) + s_{l,l}(x) = 0$ for all $x \in N$, and using \eqref{diag-plf} we have that
\begin{equation} \label{off-diag-plf-2}
is_{k,l}(x) - is_{l,k}(x) = 0 \text{ for all $x \in N$.}
\end{equation}
It follows from \eqref{off-diag-plf-1} and \eqref{off-diag-plf-2} that $s_{k,l}(x) = 0$ for all $x \in N$.  Therefore, since $A = (a_{k,l})_{k,l}  \in M_n(N)$, we have that $s(A) = \sum_{k,l} s_{k,l}(a_{k,l}) = 0$.  Hence $A \in N_n$.  Therefore $M_n(N) \subseteq N_n$.
\end{proof}

Suppose that $V$ is a matrix ordered $*$-vector space with matrix order unit $e$. As before, $N$ is a $*$-subspace of $V$, the quotient $V/N$ is a $*$-vector space in the natural way,
and $(V/N)_h = \{ v+N : v \in V_h \}$.  Furthermore, we may identify $M_n(V/N)$ with $M_n(V) / M_n(N)$, and we see that $(M_n(V) / M_n(N))_h = \{A + M_n(N) : A^* = A \}$.  Moreover, for any
$X \in M_{n,m}(\bb C)$ we have that $X^* M_n(N) X \subseteq M_m(N)$.  We also see that $(e+N)_n = e_n + M_n(N)$.

\begin{defn}
Let $(V, \{C_n \}_{n=1}^\infty, e)$ be a matrix ordered
$*$-vector space with matrix order unit $e$. Set
\begin{align*}
C_n^\textnormal{Arch} := \{ A + M_n(N) &\in M_n(V)/M_n(N) :\\ &(re_n +A) + M_n(N) \in C_n + M_n(N) \text{ for all $r > 0$} \}.
\end{align*}
and let $V_\textnormal{Arch} := (V/N, \{C_n^\textnormal{Arch} \}_{n=1}^\infty, e+N)$.
\end{defn}

\begin{prop} \label{Arch-is-what-claimed-prop}
Let $(V, \{C_n \}_{n=1}^\infty)$ be a matrix ordered $*$-vector space with
matrix order unit $e$.
Then $V_\textnormal{Arch} = (V/N, \{C_n^\textnormal{Arch} \}_{n=1}^\infty, e)$
is a matrix ordered $*$-vector space, and $e+N$ is an Archimedean matrix order unit for this space.
\end{prop}

\begin{proof}
By identifying, $M_n(V/N)$ with $M_n(V) / M_n(N)$ we may use Lemma~\ref{N-matrix-is-Nn-lem} to conclude that for any $n \in \mathbb{N}$ we have $$(M_n(V/N), C_n^\textnormal{Arch}, e_n+M_n(N)) = (M_n(V) / N_n, C_n^\textnormal{Arch}, e_n+N_n).$$  Thus we see that $(M_n(V/N), C_n^\textnormal{Arch}, e_n+M_n(N))$ is the Archimedeanization of the matrix ordered space $(M_n(V), C_n, e_n)$ (see \cite[Definition~3.15]{pt}).  Since the Archimedeanization is an AOU space, this implies that $C_n^\textnormal{Arch}$ is a cone, $C_n^\textnormal{Arch} \cap -C_n^\textnormal{Arch} = \{ 0 \}$, and $e_n+M_n(N)$ is an Archimedean order unit.

All that remains is to show that the family $\{ C_n^\textnormal{Arch} \}_{n=1}^\infty$ is compatible.  Suppose that $A \in C_n^\textnormal{Arch}$ and $X \in M_{n,m}(\bb C)$.
Since $X^* e_nX \in M_m(V)$ and $e$ is a matrix order unit, it follows that there exists $r_0 > 0$ such that $r_0 e_m - X^* e_n X \in C_m$.  Since $A \in C_n^\textnormal{Arch}$ we have that $(r e_n+A) + M_n(N) \in C_n + M_n(N)$ for all $r >0$.  Hence for any $r >0$ we have that $$\left(\frac{r}{r_0} e_n+A\right) + M_n(N) \in C_n + M_n(N)$$ and since $X^* C_nX \subseteq C_m$ and $X^*M_n(N)X \subseteq M_m(N)$ we have $$X^*\left(\frac{r}{r_0} e_n+A\right)X + M_m(N) \in C_m + M_m(N)$$ or $$\left(\frac{r}{r_0} X^*e_nX+X^*AX \right) + M_m(N) \in C_m + M_m(N).$$  By adding $re_m - (r/r_0) X^*e_nX = (r/r_0)(r_0e_m - X^*e_n X) \in C_m$ to this element we obtain $$\left(re_m +X^*AX \right) + M_m(N) \in C_m + M_m(N).$$
Since this holds for all $r >0$ we have that $X^*AX  +M_m(N) \in C_m^\textnormal{Arch}$, and it follows that $X^*C_n^\textnormal{Arch}X \subseteq C_m^\textnormal{Arch}$.  We thus conclude that $V_\textnormal{Arch} = (V/N, \{C_n^\textnormal{Arch} \}_{n=1}^\infty, e)$ is a matrix ordered $*$-vector space with Archimedean order unit $e+N$.
\end{proof}

\begin{remark}
As described in the proof of Proposition~\ref{Arch-is-what-claimed-prop}, it is useful to realize that the Archimedeanization of a matrix ordered space $(V, \{C_n \}_{n=1}^\infty, e)$ is obtained by forming the Archimedeanization of $(M_n(V), C_n, e_n)$ at each matrix level.
\end{remark}

In addition, we have the following matricial version of \cite[Theorem~3.16]{pt}.

\begin{thm} \label{complex-Arch-charact-thm}
Let $(V, \{C_n \}_{n=1}^\infty, e)$ be a matrix ordered
$*$-vector space with matrix order unit $e$, and let $V_\textnormal{Arch}$ be the
Archimedeanization of $V$ with Archimedean matrix order unit $e+N$.
Then there exists a unital surjective completely positive linear
map $q : V \to V_\textnormal{Arch}$ with the property that whenever
$(W, \{C_n'\}_{n=1}^\infty, e' )$ is an operator system with Archimedean
order unit $e'$, and $\phi : V \to W$ is a unital completely positive
linear map, then there exists a unique completely positive linear map
$\tilde{\phi} : V_\textnormal{Arch} \to W$ with $\phi = \tilde{\phi} \circ q$.
\begin{equation*}
\xymatrix{  V_\textnormal{Arch} \ar@{-->}[rd]^{\tilde{\phi}} & \\ V \ar[u]^q  \ar[r]^<>(.35)\phi & W}
\end{equation*}
Moreover, this property characterizes $V_\textnormal{Arch}$:
if $V'$ is any ordered $*$-vector space with an Archimedean order unit and
$q' : V \to V'$ is a unital surjective positive linear map with the above
property, then $V'$ is isomorphic to $V_\textnormal{Arch}$ via a
unital complete order isomorphism.
\end{thm}

\begin{proof}
Let $q : V \to V/N$ be the quotient map $q(v) = v + N$.
Let $(W, \{C_n'\}_{n=1}^\infty )$ be an operator system with
Archimedean order unit $e'$, and $\phi : V \to W$ be a unital
completely positive linear map, and choose any $v \in N$.
We see that for any state $f : W \to \bb C$ the map
$f \circ \phi : V \to \bb C$ is a state on $V$, and hence by
the definition of $N$ we have that $f(\phi(v)) = 0$.
Since $(W,C_1',e')$ is an AOU space and
$f(\phi(v)) = 0$ for all states $f : W \to \bb C$ it follows from
\cite[Proposition~3.12]{pt} that $\phi(v) = 0$.
Hence $\phi$ vanishes on $N$, and the map
$\tilde{\phi} : V/N \to W$ given by $\tilde{\phi} (v+N) = \phi(v)$
is well defined and makes the above diagram commute.
Note also that $\tilde{\phi}_n (A + M_n(N)) = \phi_n(A)$.

Furthermore, if $A+M_n(N) \in C_n^\textnormal{Arch}$, then $(re_n+A) + M_n(N) \in C_n + M_n(N)$ for all $r > 0$.
Applying $\tilde{\phi}_n$ gives that $re_n' + \phi_n(A) \in C_n'$ for all $r>0$
(recall that $\phi$ is completely positive so that $\phi_n(C_n) \subseteq C_n'$).
Since $e'$ is an Archimedean matrix order unit, this implies that $\phi_n(A) \in C_n'$ and thus $\tilde{\phi}_n (A + M_n(N)) \in C_n'$.  Hence $\tilde{\phi}$ is completely positive.

Finally, to see that $\tilde{\phi}$ is unique, simply note that any $\psi : V_\textnormal{Arch} \to W$ that makes the above diagram commute would have $\psi (v + N) = \psi ( q (v)) = \phi (v) = \tilde{\phi} (q (v)) = \tilde{\phi} (v+N)$ so that $\psi = \tilde{\phi}$.

The fact that $V_\textnormal{Arch}$ is characterized up to unital complete order isomorphism by the universal property follows from a standard diagram chase.
\end{proof}

\begin{remark} \label{Arch-special-case-rem}
The special case that is of interest to us is the case when
$(V, \{C_n\}_{n=1}^\infty, e)$ is a matrix ordered $*$-vector space with
matrix order unit $e$, and $(V, C_1, e)$ is an AOU space. In this case, since $e$ is
an Archimedean order unit for the ground level $n=1$, we have
that $N = \{0\}$, $V/N = V$, and $C_1^\textnormal{Arch} = C_1$.
Furthermore, for the higher levels $n \geq 2$, the fact that $N = \{ 0 \}$ shows that
$C_n^\textnormal{Arch} = \{A \in M_n(V) : re_n + A \in C_n \mbox{ for all } r > 0\}$.
Thus we see that in this case each $C_n^\textnormal{Arch}$ is obtained by enlarging $C_n$.
We shall show in the next proposition that each $C_n^\textnormal{Arch}$
may be viewed as the closure of $C_n$ in the order topology on $M_n(V)$.
\end{remark}

\begin{prop}
Let $(V, \{C_n \}_{n=1}^\infty)$ be a matrix ordered $*$-vector
space with matrix order unit $e$, and suppose that $(V, C_1, e)$
is an AOU space. Then $V_\textnormal{Arch}$ is the operator system
with underlying space $V$, matrix ordering
$\{C_n^\textnormal{Arch} \}_{n=1}^\infty$ given by
$C_1^\textnormal{Arch} = C_1$ and
$$C_n^\textnormal{Arch} := \{ A \in M_n(V) : re_n + A \in C_n \mbox{ for all } r > 0 \} \quad \mbox{ for } n \geq 2,$$
together with the Archimedean matrix order unit $e$.  In addition, each
$C_n^\textnormal{Arch}$ is equal to the closure of $C_n$ in the order topology on $M_n(V)$.
\end{prop}
\begin{proof}
The fact that $V_\textnormal{Arch}$ is equal to $V$ with the matrix
orderings given above follows from Remark~\ref{Arch-special-case-rem}.
Thus we need only show that $C_n^\textnormal{Arch}$ is equal to the closure
of $C_n$ in the order topology on $M_n(V)$.  Since $M_n(V)_h$ is closed in
$M_n(V)$ in the order topology, and $C_n \subseteq M_n(V)_h$, it suffices
to show that $C_n^\textnormal{Arch}$ is equal to the closure of $C_n$ in
the order topology on $M_n(V)_h$.
However, this follows from \cite[Theorem 2.34]{pt}.
\end{proof}

We now return to the case of interest.

\begin{defn}\label{dcmax}
Let $(V, V^+, e)$ be an AOU space.  We set
$$C_n^{\max}(V) = \{ A \in M_n(V) : re_n +A \in D_n^{\max}(V) \mbox{ for all } r > 0\},$$
$\cl C^{\max}(V) = \{ C_n^{\max}(V) \}_{n=1}^{\infty}$\textbf{}
and define $$\OMAX(V)= (V, \cl C^{\max}(V), e).$$
\end{defn}

The following theorem summarizes the consequences of the above results. The last statement
can be proved in the same manner as the last statement of Theorem \ref{autocp}.

\begin{thm}\label{unomax}
Let $(V,V^+,e)$ be an AOU space.
\begin{itemize}
\item[(i)] $\OMAX(V)$ is an operator system structure on $(V,V^+,e).$
\item[(ii)] If $(V, \{ P_n \}_{n=1}^{\infty}, e)$ is any operator system structure on $(V,V^+,e),$ then $C_n^{\max}(V) \subseteq P_n,$ for all $n \ge 1.$
\item[(iii)]  If $\cl S$ is any operator system and $\phi: V \to \cl S$ is a unital positive map, then $\phi: \OMAX(V) \to \cl S$ is completely positive.
\end{itemize}
Moreover, if $\tilde{V} = (V,\{C_n\}_{n=1}^{\infty},e)$
is an operator system structure on $V$ with $C_1 = V^+$ and such that for every operator system $W$
any unital positive map $\Psi : \tilde{V}\rightarrow W$ is
completely positive, then the identity map is a unital complete order
isomorphism from $\tilde{V}$ onto $\OMAX(V)$.
\end{thm}


\section{The Matricial State Space}\label{mss}

A matricial order on a $*$-vector space induces a natural matrix order
on its dual space.  In this section we describe
the correspondence
between the various operator system structures that an AOU space
can be endowed with and the corresponding matricial state spaces.
Unfortunately, duals of AOU spaces are not in general AOU spaces,
but they are normed $*$-vector spaces.
As was shown in \cite{pt}, the order norm on the self-adjoint part
$V_h$ of an AOU space $V$ has many possible extensions to a norm on $V,$ but all
these norms are equivalent and hence the set of continuous linear functionals
on $V$ with respect to any of these norms coincides with the same space
which we shall denote by $V^{\prime}.$
For a functional $f\in V'$ we let $f^*\in V'$ be the functional given by
$f^*(v) = \overline{f(v^*)}$; the mapping $f\rightarrow f^*$ turns
$V'$ into a $*$-vector space.

\begin{defn}\label{d_fij}
Given an operator system structure
$\{P_n \}_{n=1}^{\infty}$ on an AOU space $(V, V^+,e)$, set
$$P_n^d = \{f: M_n(V) \to \bb C \ : \ f \mbox{ linear and } f(P_n) \subseteq \bb R^+\}.$$
Given $f \in P_n^d,$ we define $f_{i,j} :V \to \bb C,$
by $f_{i,j}(v) = f(v \otimes E_{i,j})$.
\end{defn}

\begin{lemma}
Given an operator system structure
$\{P_n\}_{n=1}^{\infty}$ on an AOU space $(V,V^+,e)$ and $f \in P_n^d$
then the functionals $f_{i,j}$ from Definition \ref{d_fij} belong to $V^{\prime}.$
\end{lemma}
\begin{proof} Since $C_n^{\max} \subseteq P_n \subseteq C_n^{\min},$ we have that $(C_n^{\min})^d \subseteq P_n^d \subseteq (C_n^{\max})^d,$ so it
suffices to prove that $f_{k,l} \in V^{\prime},$ whenever
$f \in (C_n^{\max})^d.$

Note that $C_1^{\max} = V^+,$ so we must first show that $(V^+)^d \subseteq V^{\prime}.$ To this end, let $f: V \to \bb C$ with $f(V^+) \ge 0$.  Let $f(e) =t \ge 0$.  Then for any $v \in V_h$ with $-re \le v \le +re$ ($r > 0$), we have that
$-rt \le f(v) \le +rt.$  Hence $|f(v)| \le t \|v\|_{m},$ and it follows that
$f$ is continuous on $V_h$. Since every element $v\in V$ can be written in the
form $v = \frac{v+v^*}{2} + i\frac{v-v^*}{2i}$ and the mapping $v\rightarrow v^*$ is
continuous, we have that $f\in V'$.

Next, if $f \in (C_n^{\max})^d$ and $v \in V^+$ then
by the definition of $C_n^{\max}$ we have that
$v \otimes E_{k,k} \in C_n^{\max}.$
Hence, $f_{k,k}(v) \ge 0.$  Thus, $f_{k,k} \in (V^+)^d \subseteq V^{\prime}.$

Finally, one checks that if $v \in V^+,$ then the four elements,
$v \otimes (E_{k,k} + E_{k,l} +E_{l,k} + E_{l,l}),
v \otimes (E_{k,k} - E_{k,l} - E_{l,k} + E_{l,l}), v
\otimes (E_{k,k} + iE_{k,l} -iE_{l,k} + E_{l,l}),$ and
$v \otimes (E_{k,k} -iE_{k,l} +iE_{l,k} + E_{l,l})$ are all in $C_n^{\max},$
from which it follows that $f_{k,k} + f_{k,l} + f_{l,k} +f_{l,l},
f_{k,k} - f_{k,l} -f_{l,k}+f_{l,l}, f_{k,k} +if_{k,l} -if_{l,k} + f_{l,l},$
and $f_{k,k} -if_{k,l} +if_{l,k} + f_{l,l}$ are all in $(V^+)^d.$
Thus, $f_{k,l} \in V^{\prime},$ by taking linear combinations.
\end{proof}

Identifying each $f \in P_n^d$ with $(f_{i,j}) \in M_n(V^{\prime}),$
we shall regard $P_n^d$ as sitting inside $M_n(V^{\prime}).$ Conversely, we
identify $(f_{i,j}) \in M_n(V^{\prime})$ with the linear functional
$f:M_n(V) \to \bb C,$ given by $f((v_{i,j})) = \sum_{i,j} f_{i,j}(v_{i,j}).$

\begin{thm}
Let $\{ P_n \}_{n=1}^{\infty}$ be an operator system structure on
the AOU space $(V, V^+,e).$
Then $\{ P_n^d \}_{n=1}^{\infty}$ is a matrix ordering on
the ordered $*$-vector space $V^{\prime}$ with
$P_1^d=(V^+)^d.$ Conversely, if $\{ Q_n \}_{n=1}^{\infty},$
is any matrix ordering on the $*$-vector space $V^{\prime}$
with $Q_1= (V^+)^d$ and we set
$$^dQ_n = \{ v \in M_n(V): f(v) \ge 0 \text{ for all } f \in Q_n \},$$
then $\{ ^dQ_n \}_{n=1}^{\infty}$ is an operator system structure on $(V,V^+,e).$
\end{thm}
\begin{proof}  We leave it to the reader to check the
first claim and we focus on the proof that $\{ ^dQ_n \}$
defines an operator system structure on $(V,V^+,e).$
First, note that $^dQ_1 = V^+.$

Let us denote by $\overline{Y}$ the transpose of the adjoint
$Y^*$ of a matrix $Y\in M_{n,m}$.
It is easy to see that each $^dQ_n$ is a cone in $M_n(V).$
We now show that $^dQ_n \cap (-^dQ_n) = \{0\}.$ Let $v=(v_{i,j})$
be in the intersection and let $s \in S(V) \subseteq (V^+)^d.$
Let $X$ be the column
vector with entries $\alpha_1,\dots,\alpha_n \in \bb C$ (in this order).
Since $XQ_1X^* \subseteq Q_n,$ we
have that  $f = (\alpha_i s \overline{\alpha_j}) \in Q_n.$
Hence, $0 = f(v) = \sum_{i,j} \alpha_i s(v_{i,j}) \overline{\alpha_j}.$
But this latter quantity is equal to
$((s(v_{i,j}))\overline{X},\overline{X})$.
Since $X$ was an arbitrary vector, we have that
$(s(v_{i,j})) =0,$ and since $s$ was arbitrary and the states
separate points in $V$ \cite[Proposition~3.12]{pt}, we have that $(v_{i,j}) =0.$
Thus $^dQ_n \cap (-{}^dQ_n) = \{0\}.$

Fix $X \in M_{n,m}(\bb C)$. We show that $X^*(^dQ_n)X \subseteq$  $^dQ_m.$
To this end, observe that if $v=(v_{i,j}) \in$ $^dQ_n$
and $f=(f_{k,l}) \in Q_m,$ then $f(X^*vX) = (\overline{X}f\overline{X}^*)(v) \ge 0,$
since $\overline{X} f \overline{X}^* \in Q_n$ and the family
$\{Q_n\}_{n=1}^{\infty}$ is a matrix ordering.
Thus, we have shown that the family $\{^dQ_n\}_{n=1}^{\infty}$
is a matrix ordering.
It follows from Proposition \ref{minmax} (iii) that
$e$ is a matrix order unit for
$\{^dQ_n\}_{n=1}^{\infty}$.

Finally, to see that $e$ is an Archimedean matrix order unit, observe that if $v \in M_n(V)$ and $re_n +v \in$ $^dQ_n$ for all $r>0,$ then $rf(e_n) + f(v) \ge 0$ for all $r>0$ and all $f \in Q_n$. Hence $f(v) \ge 0$ for all $f \in Q_n,$ and it follows that $v \in$ $^dQ_n$.
\end{proof}

Note that the weak*-topology on $V^{\prime}$ endows $M_n(V^{\prime})$ with a topology
which coincides with the weak*-topology that comes from the identification of
$M_n(V^{\prime})$ with the dual of $M_n(V).$ Thus, we shall refer to this topology,
unambiguously, as {\it the weak*-topology on $M_n(V^{\prime}).$}

\begin{thm}\label{th_donet}
Let $(V,V^+,e)$ be an AOU space.
The mappings $P_n \to P_n^d$ and $Q_n \to {}^dQ_n$
establish a
one-to-one inclusion-reversing correspondence between operator system structures
$\{ P_n \}_{n=1}^{\infty}$ on $(V,V^+,e)$ and matrix orderings
$\{ Q_n \}_{n=1}^{\infty}$ on $V^{\prime}$ with $Q_1 = (V^+)^d$
for which each $Q_n$ is weak*-closed.
\end{thm}
\begin{proof} It suffices to show that $^d(P_n^d) = P_n$ and $(^dQ_n)^d = Q_n,$ which follows from a standard dual cone argument.
\end{proof}

\begin{remark}\label{r_wst}
In Theorem~\ref{w*-closure-Cmin} we identify the dual cones of the minimal and maximal matrix orderings.
Although the cone $Q_n^{\min}$ defined below
is not weak*-closed, we show in Theorem~\ref{w*-closure-Cmin} that
$^dQ_n^{\min} = C_n^{\min}$, from which it follows that
$(C_n^{\min})^d$ is the weak*-closure of $Q_n^{\min}$.
\end{remark}

\begin{defn} Let $(V,V^+,e)$ be an AOU space. Set
$$Q_n^{\min} = \{X^*GX : X\in M_{m,n}, G = {\rm diag}(g_1,\dots,g_m), g_i \in (V^+)^d, m\in\bb{N}\}$$
and
$$Q_n^{\max} = \{ (f_{i,j}) \in M_n(V^{\prime})
: (f_{i,j}(v)) \in M_n^+ \text{ for all } v \in V^+ \}.$$
\end{defn}

\begin{remark}\label{altqnmin} Using the same technique as the one used in the
  proof of Proposition~\ref{minmax}(ii), one easily shows that
$$Q_n^{\min} = \{ \sum_{i=1}^k P_i \otimes g_i : g_i \in (V^+)^d, P_i
\in M_n^+, i=1,\dots,k, k \in \bb N \}.$$

\end{remark}

\begin{thm} \label{w*-closure-Cmin}
Let $(V,V^+,e)$ be an AOU space. Then $^dQ_n^{\min} = C_n^{\min}$ and $(C_n^{\max})^d = Q_n^{\max}.$
\end{thm}
\begin{proof}
It is easily checked that both
$\{ Q_n^{\min} \}$ and $\{ Q_n^{\max} \}$ are
matrix orderings on $V^{\prime}$ with $Q_1^{\min} = Q_1^{\max} = (V^+)^d.$

Let $v=(v_{i,j}) \in C_n^{\min},$ and let $f=X^*GX \in Q_n^{\min}$.
Then $f(v) = G(\overline{X}v \overline{X}^*) = \sum_i g_i(w_{i,i}),$
where $(w_{i,j}) = \overline{X}v\overline{X}^* \in C_m^{\min}.$
But this sum is non-negative since $w_{i,i} \in V^+$ for each $i$.
Thus $C_n^{\min} \subseteq$ $^dQ_n^{\min},$ and equality follows from the fact
that they are both operator system structures on
$(V,V^+,e)$ and $C_n^{\min}$ is the largest possible.

Let $f=(f_{i,j}) \in (C_n^{\max})^d,$ let $v \in V^+,$ and
let $X= ( \alpha_1,\dots,\alpha_n) \in M_{1,n}(\bb C).$
Then $\langle (f_{i,j}(v))X^*,X^* \rangle =
\sum_{i,j} f_{i,j}(v) \alpha_i \overline{\alpha_j} = f(XvX^*) \ge 0.$
Since $X$ was arbitrary, this shows that $(f_{i,j}(v)) \in M_n^+,$
and since $v$ was arbitrary, $(f_{i,j}) \in Q_n^{\max}.$
Thus $(C_n^{\max})^d \subseteq Q_n^{\max}.$
But since $\{ ^dQ_n^{\max} \}$ is an operator system structure,
we have that $C_n^{\max} \subseteq$ $^dQ_n^{\max},$ and hence
$(^dQ_n^{\max})^d \subseteq (C_n^{\max})^d$.
On the other hand, $Q_n^{\max}$
is easily seen to be weak*-closed.  Hence
$(^dQ_n^{\max})^d = Q_n^{\max}$, and equality follows.
\end{proof}


\section{Comparison of various structures}\label{s_cvs}

Given a unital $C^*$-algebra or, more generally,
an operator system $\cl S$, at ground level it is an AOU space,
so we may form new operator systems, $\OMIN(\cl S)$ and $\OMAX(\cl S).$
Also, since it is a normed space, we may form the operator spaces
$\MIN(\cl S)$ and $\MAX(\cl S)$.  In this section, we compare
these structures, describing when they are identical and, more
generally, when the identity map between these various structures is a
completely bounded isomorphism.  We have already seen one result of this
type: Proposition~\ref{p_minomin} states that, for every AOU space $(V,V^+,e)$, the
identity map from $\MIN(V_{\min})$ to $\OMIN(V)$ is a complete isometry.
We begin by identifying the norm structure on $\OMAX(V).$

In \cite[\S4]{pt} three norms on AOU spaces that extend the order norm $\| \cdot \|$ on the self-adjoint elements were studied:
the {\em minimal norm} $\| \cdot\|_m$, the {\em maximal norm}
$\| \cdot \|_M,$ and the {\em decomposition norm} $\| \cdot \|_{\dec}.$
We have encountered the minimal and maximal norms earlier. The
decomposition norm is given by the following formula: for $v \in V,$
we set
$$\|v\|_{\dec} = \inf \left\{ \left\| \sum_{j=1}^k |\lambda_j| v_j \right\| : v = \sum_{j=1}^k
\lambda_j v_j \text{ with } \lambda_j \in \bb C \text{ and } v_j \in V^+ \right\}.$$
The following result shows that the decomposition norm
is the largest norm to come from an operator system structure on $V.$

\begin{prop}
Let $(V,V^+,e)$ be an AOU space.
Then for $v \in V,$ we have that
\begin{equation}\label{seq}
\|v\|_{\dec} = \|v\|_{\OMAX(V)}= \sup \{ \|\phi(v) \| \ : \ \phi:V \to B(H) \},
\end{equation}
where the supremum is taken over all Hilbert spaces and
all unital positive maps $\phi.$
\end{prop}
\begin{proof}
Suppose that $\phi : V\rightarrow B(H)$ is a unital positive
map. By Theorem \ref{unomax}, $\phi$ is completely positive
as a map from $\OMAX(V)$ into $B(H)$ and hence it is completely
contractive. It follows that $\phi$ is contractive and hence
$\|\phi(v)\| \leq \|v\|_{\OMAX(V)}$, for each $v\in V$. On the other hand,
if $\phi : \OMAX(V)\rightarrow B(H)$ is a unital complete isometry
then $\phi$ is completely positive and $\|v\|_{\OMAX(V)} = \|\phi(v)\|$ for
each $v\in V$. The second equality in (\ref{seq}) follows.

Let $\phi : V\rightarrow B(H)$ be a unital positive map. Fix $v\in V$
and write $v = \sum_{i=1}^k \lambda_i v_i$ where $v_i\in V^+$ and
$\lambda_i\in\bb{C}$, $i = 1,\dots,k$.
Using \cite[Proposition~5.10]{pt} and Theorem \ref{unomax} we have
\begin{align*}
\|\phi(v)\|_{B(H)}&  =
\left\|\sum_{i=1}^k \lambda_i \phi(v_i)\right\|_{B(H)} \leq
\left\|\sum_{i=1}^k |\lambda_i| \phi(v_i)\right\|_{B(H)}\\
& = \left\|\phi\left(\sum_{i=1}^k |\lambda_i|
v_i\right)\right\|_{B(H)}
\leq \left\|\sum_{i=1}^k |\lambda_i| v_i\right\|_{\OMAX(V)}
= \left\|\sum_{i=1}^k |\lambda_i| v_i\right\|.
\end{align*}
Taking the infimum over all representations of $v$ of the form $v = \sum_{i=1}^k \lambda_i v_i$,
$v_i\in V^+$, and supremum over all unital positive maps $\phi : V\rightarrow B(H)$,
we obtain that
$\|v\|_{\OMAX(V)}\leq \|v\|_{\dec}$.

It remains to establish the inequality $\|v\|_{\dec}\leq \|v\|_{\OMAX(V)}$.
Suppose that $v\in V$ is such that $\|v\|_{\OMAX(V)}\leq 1$.
Then $\left(\smallmatrix e & v\\ v^* & e\endsmallmatrix\right) \in C_2^{\max}$.
Consider first the special case where the matrix
$\left(\smallmatrix e & v\\ v^* & e\endsmallmatrix\right)$
lies in the smaller cone $D_2^{\max}$. Then there exist positive matrices
$a_j = \left(\smallmatrix \alpha_j & \lambda_j\\ \overline{\lambda_j} & \beta_j\endsmallmatrix\right)\in M_2^+$
and elements $v_j\in V^+$, $j = 1,\dots,k$, such that
$\left(\smallmatrix e & v\\ v^* & e\endsmallmatrix\right) = \sum_{j=1}^k a_j\otimes v_j$.
This implies that
$$e = \sum_{j=1}^k \alpha_j v_j = \sum_{j=1}^k \beta_j v_j, \ \ v = \sum_{j=1}^k \lambda_j v_j.$$
We will show that
\begin{equation}\label{inmod}
\left\|\sum_{j=1}^k |\lambda_j| v_j\right\| \leq 1.
\end{equation}

Conjugating each $a_j$ by a diagonal unitary, it is easily seen that
$b_j = \left(\smallmatrix \alpha_j & |\lambda_j| \\ |\lambda_j| &
  \beta_j \endsmallmatrix \right) \in M_2^+.$  Hence,
$$A= \begin{pmatrix} e & \sum_{j=1}^k |\lambda_j|v_j \\ \sum_{j=1}^k
  |\lambda_j|v_j & e \end{pmatrix} = \sum_{j=1}^k b_j \otimes v_j \in
D_2^{\max}.$$
If $\phi: V \to \bb C$ is any state, then $\phi$ is completely
positive on $\OMAX(V),$ and hence $\phi_2(A) \in M_2^+$ from which
it follows that $| \phi(\sum_{j=1}^k |\lambda_j| v_j) | \le 1.$
Taking the supremum over all states $\phi$, we obtain
$\| \sum_{j=1}^k |\lambda_j| v_j
\| \le 1$, and consequently  $\|v\|_{\dec} \le 1$.

Now consider the general case
$\left(\smallmatrix e & v\\ v^* & e\endsmallmatrix\right) \in C_2^{\max}$.
By Definition \ref{dcmax}, we have that
$\left(\smallmatrix (r+1)e & v\\ v^* & (r+1)e\endsmallmatrix\right) \in D_2^{\max}$
for each $r > 0$.
Thus,
$\left(\smallmatrix e & \frac{1}{r+1}v\\ \frac{1}{r+1}v^* & e\endsmallmatrix\right) \in D_2^{\max}$
for each $r > 0$. By the previous paragraph,
$\left\|\frac{1}{r+1}v\right\|_{\dec} \leq 1$ for each $r > 0$ which implies that
$\|v\|_{\dec} \leq 1$.
\end{proof}

\begin{prop}\label{cxominv}
Let $X$ be a compact Hausdorff space and let $C(X)$ be the
$C^*$-algebra of continuous functions on $X.$
Then the identity map is a complete order
isomorphism between $C(X), \OMIN(C(X))$ and $\OMAX(C(X)).$
\end{prop}
\begin{proof} The identification of $C(X)$ with
$\OMIN(C(X))$ follows from the fact that every unital positive map
from an operator system into $C(X)$ is completely positive \cite[Theorem~3.9]{Pa}
and the characterization
of $\OMIN(C(X))$ given in Theorem \ref{autocp}.
The identification of $C(X)$
with $\OMAX(C(X))$ follows from the fact that every unital positive map
from $C(X)$ into an operator system is completely positive
\cite[Theorem~3.11]{Pa} and the characterization of $\OMAX(C(X))$ given in Theorem~\ref{unomax}.
\end{proof}

\begin{prop}
Let $\cl S$ be an operator system. Then the identity map from
$\OMIN(\cl S)$ to $\cl S$ is completely bounded with $\|\id\|_{cb} =C$
if and only if for every operator system $\cl T,$ every unital
positive map from $\cl T$ into $\cl S$ is completely bounded
and the supremum of the completely bounded norms of all such maps is $C.$
\end{prop}

\begin{proof} Assume that the map
$\id: \OMIN(\cl S) \to \cl S$ is completely bounded.
Given a unital positive map
$\phi: \cl T \to \cl S,$ we may write
$\phi = \id \circ \gamma,$ where $\gamma$ is the same map as
$\phi$ but with $\OMIN(\cl S)$ regarded as its range.
Since $\gamma$
is unital and positive, by Theorem \ref{autocp},
it is unital and completely positive,
and hence it is completely contractive.
Thus, $\phi$ is the composition of
two completely bounded maps and hence is completely bounded with
$\|\phi\|_{cb} \le \|\id\|_{cb} = C.$

The converse implication follows by taking
$\cl T = \OMIN(\cl S)$ and $\phi = \id$.
\end{proof}

\begin{prop}
Let $\cl S$ be an operator system.
Then the following are equivalent:
\begin{itemize}
\item[(i)] the identity map $\id$ from $\cl S$ to
$\OMAX(\cl S)$ is completely bounded with $\|\id\|_{cb} =C$;
\item[(ii)] for every operator system $\cl T,$ every unital positive map
from $\cl S$ into $\cl T$ is completely bounded and the supremum
of the completely bounded norms of all such maps is $C$;
\item[(iii)] for every Hilbert space $H,$ every unital positive map
$\phi$ from $\cl S$ into $B(H)$ decomposes as a difference of
two completely positive maps $\phi =\phi_1 - \phi_2,$ and
$$C = \sup{}_{\phi} \ \inf\{\|\phi_1(e) + \phi_2(e)\| : \phi = \phi_1 - \phi_2\}$$
where the supremum is taken over all unital positive maps $\phi:
\cl S\rightarrow B(H)$ and the infimum over all decompositions
of $\phi$ with $\phi_1$ and $\phi_2$ completely positive.
\end{itemize}
\end{prop}
\begin{proof} The equivalence of (i) and (ii) is similar to the last proof.
The equivalence of (ii) and (iii) follows by first observing that the arbitrary
operator system $\cl T$ in (ii) can be replaced by $B(H)$ for some $H,$
since every $\cl T$ embeds into some $B(H).$  Then by Wittstock's
decomposition theorem \cite[Theorem~8.5]{Pa}, a unital positive map is
completely bounded if and only if it decomposes as a difference
of completely positive maps and the infimum of $\|\phi_1(e) + \phi_2(e) \|$
over all such decompositions is equal to $\|\phi\|_{cb}.$
\end{proof}

\begin{prop}
Let $X$ be a compact Hausdorff space and let
$\cl S \subseteq C(X)$ be an operator
system of codimension $n < \infty.$ Then the identity map from
$\cl S$ to $\OMAX(\cl S)$ is completely bounded with
completely bounded norm at most $2n+1.$
\end{prop}
\begin{proof} By \cite[Lemma~11.9]{Pa} (see also \cite{dp}) for every
  $\epsilon > 0,$ there exists a completely positive map $\phi:C(X)
  \to \cl S$ with $\|\phi\|_{cb} \le n+1+ \epsilon$ and a positive
  linear functional $s: C(X) \to \bb C$ with $\|s\| \le n+ \epsilon$ such
  that $x = \phi(x) - s(x)1,$ for every $x \in \cl S.$  Now if $\cl T$
  is an operator system and $\psi:
  \cl S \to \cl T$ is a
unital positive map, then $\psi(x) = \psi \circ \phi(x) - s(x) 1.$  Hence
$\|\psi\|_{cb} \le \|\psi \circ \phi \|_{cb} + \|s\|.$  Since $\psi
\circ \phi:C(X) \to \cl T$ is positive, it is completely positive and
$\|\psi \circ \phi\|_{cb} = \|\psi \circ \phi (1) \| \le \|\phi(1)\|
\le n+1+ \epsilon.$ Thus we have that $\psi$ has completely bounded
norm at most $2n+1+2 \epsilon$. The proof is completed by letting $\cl
T = {\rm OMAX}(\cl S)$ and letting $\psi$ be the identity map.
\end{proof}

If $W$ is a normed space, the identity map from $\MIN(W)$ to
$\MAX(W)$ is completely bounded if and only if $W$ is finite
dimensional \cite{Pa2}. However, by the above results we see
that there are plenty of infinite dimensional operator systems
for which $\MIN(\cl S) = \OMIN(\cl S)$ and the identity map from
$\OMIN(\cl S)$ to $\OMAX(\cl S)$ is completely bounded. Thus, in
general, the identity map from $\OMAX(\cl S)$ to $\MAX(\cl S)$
will not be completely bounded.
In fact, we have the following characterization of when this happens:

\begin{prop} Let $(V,V^+,e)$ be an AOU space. Then
the identity map from $\OMAX(V)$ to $\MAX(V_{\min})$ is completely bounded if and only if
for every Hilbert space $H,$ every bounded map $\phi : V_{\min} \to B(H)$
decomposes as $$\phi = (\phi_1 - \phi_2) + i(\phi_3 - \phi_4)$$
where each $\phi_j : V \to B(H)$ is positive.
\end{prop}
\begin{proof}
Assume that the decomposition of bounded maps holds,
and suppose that $\MAX(V_{\min}) \subseteq$ $B(H)$
completely isometrically for some Hilbert space $H$. Let $\phi:V_{\min} \to \MAX(V_{\min})$
be the identity map and let $\phi_j$ be a positive mapping for $j=1,2,3,4$, and
such that $\phi = (\phi_1 - \phi_2) + i(\phi_3 - \phi_4)$.
Since each $\phi_j$ is positive, $\phi_j: \OMAX(V) \to B(H)$ is completely positive and hence completely
bounded for $j = 1,2,3,4$. Hence $\phi:\OMAX(V) \to B(H)$ is completely bounded.

Conversely, if the identity map from $\OMAX(V)$ to $\MAX(V_{\min})$ is completely bounded and
$\phi:V_{\min} \to B(H)$ is bounded, then $\phi: \MAX(V_{\min}) \to B(H)$ is completely bounded
and hence $\phi:\OMAX(V) \to B(H)$ is completely bounded. Applying Wittstock's decomposition
theorem, we have that $\phi = (\phi_1 - \phi_2) + i (\phi_3 - \phi_4),$ where each
$\phi_j : \OMAX(V) \to B(H)$ is completely positive and hence positive as a map from $V$ into $B(H)$.
\end{proof}

\section{Entanglement breaking maps}\label{s_ent}

In Quantum Information Theory entangled states and
separable states are important objects of
study and there is interest in maps which are ``entanglement
breaking'' (see \cite{Hol}, \cite{HS} and \cite{Ar}). In the present section we review these concepts, relate them to our
constructions of minimal and maximal operator system structures, and discuss some generalizations.

If $n\in\bb{N}$, we will be interested in three operator system structures on $M_n$: the minimal operator system structure $\OMIN(M_n)$, the maximal operator system structure $\OMAX(M_n)$, and the
operator system structure, simply denoted $M_n$, arising from the
identification of $M_n$ with $B(\bb{C}^n)$.
The cone of positive elements of $M_n$ for any of these
operator system structures coincides with the set of all
positive definite matrices in $M_n$.

A state $s : M_n \otimes M_m \to \bb C$ is called {\bf separable} if
it is a convex combination of tensor states; i.e., if
there exist $l\in \bb{N}$, states $s_i:M_n \to \bb C$, states $t_i:M_m \to \bb C$, and
real numbers $r_i >0$ for $i = 1,\dots,l$ with $\sum_{i=1}^l r_i =1$ and such that
$s = \sum_{i=1}^l r_i s_i \otimes t_i$.
Note that if $s$ is a
state and $f_i:M_n \to \bb C$ and $g_i:M_m \to \bb C$ are positive
linear functionals for $i = 1,\dots,l$ with $s = \sum f_i \otimes g_i,$ and
if we set $r_i = f_i(I_n)g_i(I_m), s_i = (f_i(I_n))^{-1} f_i$, and $t_i = (g_i(I_m))^{-1}
g_i,$ then $\sum_{i=1}^l r_i = 1$, each $s_i$ and each $t_i$ is a state, and $s= \sum_{i=1}^l r_i s_i
\otimes t_i.$ Thus a state is separable if and only if it is a sum
of tensors of positive linear functionals.

More generally, we shall call a positive linear functional (p.l.f.)
$f : M_n \otimes M_m \to \bb C$ {\bf separable} if it is a sum of
tensors of positive linear functionals. The previous paragraph shows that
these two definitions agree for states.
We recall that,
given a completely positive map
$\phi : M_k \to M_m,$ we denote by $\phi_n: M_n \otimes M_k
\to M_n \otimes M_m$ the map given by $\phi_n((v_{i,j})) = (\phi(v_{i,j}))$.
If $s : M_n \otimes M_m\rightarrow \bb{C}$ is a p.l.f., then $s \circ
\phi_n : M_n \otimes M_k \to \bb C$ is a p.l.f.
If $s$ is a state and $\phi$ is unital, then $s \circ
\phi_n$ is also a state.
The linear map $\phi : M_k \to M_m$ is called {\bf entanglement breaking} if $s \circ
\phi_n$ is a separable state for every state $s : M_n\otimes M_m \to \bb C$ and
every $n\in \bb{N}.$

In this section we relate entanglement breaking to the minimal and
maximal operator system structures studied in the previous sections.
We also prove a duality result that explains some of the identifications that occur in this
theory. We begin with a characterization of separable states.

We identify $M_n \otimes M_m$ with $M_n(M_m).$
Let $f: M_n \otimes M_m \to \bb C$ and write
  $f=(f_{i,j}),$ where $f_{i,j}: M_m \to \bb C$ are determined by the
  relation $f((B_{i,j})) =
  \sum_{i,j} f_{i,j}(B_{i,j}).$

\begin{prop}\label{p_senm}
Let
$f: M_n \otimes M_m \to \bb C$ be a positive linear functional.
Then $f$ is separable if and only if
  $f : M_n(\OMIN(M_m)) \to \bb C$ is positive.
\end{prop}
\begin{proof}
Without loss of generality, we may assume that $f$ is a state.
Using the notation of Sections~\ref{osys} and \ref{mss}
and Remark~\ref{r_wst}, we have that $f : M_n(\OMIN(M_m)) \to \bb C$ is positive if and only if
$f\in (C_n^{\min}(M_m))^{d}$, which occurs if and only if
$f\in \overline{Q_n^{\min}(M_m)}^{w^*}$.
By Remark \ref{r_two}, this happens precisely when
$f$ is the weak*-limit of separable positive linear functionals.
By dividing each of those functionals by its norm,
we may assume that each of them is a state, and hence that
$f$ is the weak*-limit of separable states.
Since the separable states are the convex hull of a compact set, they are
also a compact set by Caratheodory's theorem.
\end{proof}

We now turn our attention to a duality result.  Recall that the dual
of a matrix ordered space is again a matrix ordered space.
Let $\delta_{i,j}: M_n \to \bb C$ be the linear functional satisfying
$$\delta_{i,j}(E_{k,l}) = \begin{cases} 1 & \text{ when } (i,j) =
(k,l)
  \\ 0 & \text{ when } (i,j) \ne (k,l) \end{cases}$$
and let $\gamma_n: M_n \to M_n^{\prime}$ be the linear isomorphism defined by
$\gamma_n(E_{i,j}) = \delta_{i,j}.$  The next result is certainly
in some sense known, but the formal statement will be useful for us in the sequel.

\begin{thm}\label{gamma}
The map $\gamma_n: M_n \to M_n^{\prime}$ is a complete order
  isomorphism of matrix ordered spaces. Consequently,
  $(M_n^{\prime}, (M_n^{\prime})^+, \tr)$ is an AOU space that is
  order isomorphic to $(M_n, M_n^+, I_n),$ where $I_n$ denotes the
  identity matrix.
\end{thm}
\begin{proof} Let $A_{i,j} \in M_k, i,j = 1,\dots,n.$  We must prove
  that $\sum_{i,j} A_{i,j} \otimes E_{i,j} \in M_k(M_n)^+$ if and only
  if $F= \sum_{i,j} A_{i,j} \otimes \delta_{i,j} \in M_k(M_n^{\prime})^+.$

Writing $F= (f_{r,s})_{r,s=1}^k,$ with each $f_{r,s} \in
  M_n^{\prime},$ we have that $F \in M_k(M_n^{\prime})^+$ if and only
  if the map $F: M_n \to M_k$ given by $F(B) = (f_{r,s}(B))$ is
  completely positive.  By a theorem of Choi (see \cite[Theorem 3.14]{Pa}),
  $F$ is completely
  positive if and only if $(F(E_{i,j})) \in M_n(M_k)^+.$
However, $(F(E_{i,j})) = \sum_{i,j} A_{i,j} \otimes E_{i,j}$
which shows that $\gamma_n$ is a complete order isomorphism.

Since $\gamma_n$ is an order isomorphism, we have that
$(M_n^{\prime}, (M_n^{\prime})^+, \gamma(I_n))$ will be an AOU space
order isomorphic to $(M_n, M_n^+, I_n).$ Finally, note that
$\gamma_n(I_n) = \sum_i \delta_{i,i} = \tr.$
\end{proof}

The AOU space $(M_n, M_n^+, I_n)$ gives rise to the minimal and
the maximal operator systems
    $\OMIN(M_n)$ and $ \OMAX(M_n)$
  and their matrix ordered dual
  spaces, $\OMIN(M_n)^{\prime}$ and $\OMAX(M_n)^{\prime}$.
On the other hand, Theorem \ref{gamma} allows one to form the operator systems $\OMIN(M_n^{\prime})$ and $\OMAX(M_n^{\prime})$ corresponding to the AOU space
$(M_n^{\prime}, (M_n^{\prime})^+, \tr)$.
The following proposition explains the relationships between these objects.

\begin{prop}\label{p_dualt}
The complete order isomorphism  $\gamma_n:M_n \to M_n^{\prime}$
gives rise to the identifications
$\OMIN(M_n)^{\prime} = \OMAX(M_n^{\prime}) = \OMAX(M_n)$ and
$\OMAX(M_n)^{\prime} = \OMIN(M_n^{\prime}) = \OMIN(M_n).$
\end{prop}
\begin{proof} For clarity, let $V = M_n^{\prime}$ and observe that,
by Proposition \ref{minmax} (ii),
  $Q_k^{\min}(M_n) = D_k^{\max}(V)$.
  As in the proof of Proposition \ref{p_senm}, we see that
the unit ball of $D_k^{\max}(V)$ is compact and by Krein-\u Smulian Theorem,
$D_k^{\max}(V)$ is closed. Thus, $D_k^{\max}(V) = C_k^{\max}(V).$ Thus, we have that
  $Q_k^{\min}(M_n) = C_k^{\max}(V),$ and so $M_k(\OMIN(M_n)^{\prime})^+ =
  M_k(\OMAX(M_n'))^+,$ and so the identity map on $M_n^{\prime}$ yields a complete order
  isometry between the matrix ordered space $\OMIN(M_n)^{\prime}$ and
  the operator system $\OMAX(M_n^{\prime}).$ Now the complete order
  isomorphism $\gamma$ allows for the identification
  $\OMAX(M_n^{\prime}) = \OMAX(M_n).$

The proof of the rest of the statement is similar.
\end{proof}

We can now prove the following.

\begin{thm}\label{th_entb}
Let $\phi : M_k \to M_m$ be a linear map.
Then $\phi$ is entanglement breaking if and only if $\phi: \OMIN(M_k)
\to M_m$ is completely positive. Moreover, if $\widetilde{M_k} = (M_k,
\{P_n\}_{n=1}^{\infty}, I)$ is any operator system structure on $M_k$
with the property that for every $m,$ a map $\phi: \widetilde{M_k} \to
M_m$ is completely positive if and only if it is entanglement
breaking, then the identity map is a complete order isomorphism from
$\widetilde{M_k}$ to $\OMIN(M_k).$
\end{thm}
\begin{proof} If $\phi: \OMIN(M_k) \to M_m$ is completely positive, then $\phi^{\prime}:
  M_m^{\prime} \to \OMIN(M_k)^{\prime}$ is completely positive. Thus
  if $f=(f_{i,j}) \in M_n(M_m^{\prime})^+$, then $(\phi'(f_{i,j})) \in
  M_n(\OMIN(M_k)^{\prime})^+.$ By Proposition~\ref{p_dualt} and the definition
  of $\OMAX(M_k^{\prime})$, the state $(\phi'(f_{i,j}))$ on $M_n \otimes M_k$ is separable.
  Thus, for every
  p.l.f. $f$ on $M_n \otimes M_m,$ we have that $f \circ \phi_n =
  (\phi'(f_{i,j}))$ is a separable state.

Conversely, suppose that $\phi$ is entanglement breaking. Then
$(\phi(f_{i,j}))$ belongs to $M_n(\OMIN(M_k)^{\prime})^+$ for every
$(f_{i,j}) \in M_n(M_m^{\prime})^+$, and hence the map $\phi^{\prime}:
M_m^{\prime} \to \OMIN(M_k)^{\prime}$ is completely positive.  It
follows that $\phi: \OMIN(M_k) \to M_m$ is completely positive.

Finally, the last statement is equivalent to the assertion that the
set of completely positive maps from $\OMIN(M_k)$ into $M_m$ and from
$\widetilde{M_k}$ into $M_m$ coincides for every $m.$ This is
equivalent to equality of the dual cones, i.e., that
$C_m^{\min}(M_k)^d =  P_m^d,$ for all $m.$ By Theorem~\ref{th_donet},
the equality of the dual cones implies equality of the cones, i.e.,
$C_m^{\min}(M_k) = P_m,$ for all $m,$ and hence the identity map is a
complete order isomorphism. \end{proof}

We can use the maps $\gamma_n$ to identify the adjoint of a map
between matrix algebras as a map between matrix algebras.

\begin{defn} For a linear map $\phi:M_k \to M_m,$ we set $\phi^{\flat} = \gamma_k^{-1} \circ \phi^{\prime} \circ \gamma_m: M_m \to M_k.$
\end{defn}

We note that the map $\phi^{\flat}$ differs from the map $\phi^{\dag}$
encountered in Quantum Information Theory.
The difference has to do with the fact that our
identification of $M_k$ and $M_k^{\prime}$ is linear instead of
conjugate linear. If $A \in M_{m,k},$ $B \in M_{k,m}$ and we define
$\phi : M_k \to M_m$ by $\phi(X) = AXB,$ then it is well known that
$\phi^{\dag}(Y) = A^*YB^*.$  We will show that $\phi^{\flat}(Y) =
A^tYB^t.$ First observe that if $f : M_k \to \bb C$ is written as $f =
\sum_{i,j} y_{i,j} \delta_{i,j}$ so that $\gamma_k^{-1}(f) = Y=
(y_{i,j}) \in M_k,$  then for $X= (x_{i,j}) \in M_k,$ we have that
$f(X) = \sum_{i,j} y_{i,j}x_{i,j} = \tr(XY^t).$

Write $\phi =
(f_{i,j}),$ where $f_{i,j}: M_k \to \bb C$, $i,j = 1,\dots,m$.
We have that
$\phi^{\flat}(E_{i,j}) = \gamma_k^{-1} \circ
\phi^{\prime}(\delta_{i,j}) = \gamma_k^{-1}(f_{i,j}).$ Now
$f_{i,j}(X) = \tr(\phi(X) E_{j,i}) = \tr(XBE_{j,i}A)=
\tr(X(A^tE_{i,j}B^t)^t),$ and so $\phi^{\flat}(E_{i,j}) =
A^tE_{i,j}B^t,$ and the claim follows from the previous paragraph.

\begin{cor}\label{c_mflat}
Let $\phi: M_k \to M_m$ be a linear map.
Then $\phi: M_k \to \OMAX(M_m)$ is
  completely positive if and only if $\phi^{\flat}: M_m \to M_k$ is
  entanglement breaking. Moreover, if $\widetilde{M_m} = (M_m, \{P_n
  \}_{n=1}^{\infty}, I)$ is any operator system structure on $M_m$
  with the property that for every $k,$ a map $\phi:M_k \to
  \widetilde{M_m}$ is completely positive if and only if $\phi^{\flat}$ is
  entanglement breaking, then the identity map is a complete order
  isomorphism from $\widetilde{M_m}$ to $\OMAX(M_m).$
\end{cor}
\begin{proof} We have that $\phi: M_k \to \OMAX(M_m)$ is completely
positive if and only if $\phi^{\prime}: \OMAX(M_m)^{\prime} \to
M_k^{\prime}$ is completely positive.
Using the identifications of Proposition~\ref{p_dualt}, we see that this
happens if and only if $\phi^{\flat}:\OMIN(M_m) \to M_k$ is
completely positive, which by Theorem~\ref{th_entb} happens precisely when
$\phi^{\flat}$ is entanglement breaking.

The last statement is equivalent to the assertion that for every $k,$
the completely positive maps from $M_k$ into $\widetilde{M_m}$ and
from $M_k$ into $\OMAX(M_m)$ coincide.  Recall that a map $\phi: M_k \to \cl S,$
where $\cl S$ is an operator system, is completely positive if and
only if $(\phi(E_{i,j})) \in M_k(\cl S)^+.$ Thus, the equality of these
sets of completely positive maps, ensures that $M_k(\OMAX(M_m))^+ =
P_k$ for every $k,$ i.e., that the identity map is a complete order isomorphism.
\end{proof}

We shall call a map $\phi:M_k \to M_m$ such that $\phi^{\flat}$ is
entanglement breaking a {\bf co-entanglement breaking map.}

\begin{remark} {\rm In Theorem~\ref{entbr} we will prove that a map is
    entanglement breaking if and only if it is co-entanglement
    breaking. Consequently, in the last statement of
    Corollary~\ref{c_mflat} the map $\phi^{\flat}$ can be replaced by
    $\phi.$}
\end{remark}

\begin{prop}\label{pp}
Let $\phi:M_k \to M_m$ be a linear map.
Then $\phi:M_k \to \OMAX(M_m)$ is completely positive if and only if
there exist positive linear functionals
$s_l: M_k \to \bb C$ and matrices
$P_l \in M_m^+,$  $l = 1,\dots,q$, such that $\phi(X) = \sum_{l=1}^q s_l(X)P_l.$
\end{prop}
\begin{proof} We have that $\phi:M_k \to \OMAX(M_m)$ is completely positive
if and only if $(\phi(E_{i,j})) \in M_k(\OMAX(M_m))^+ =
C_k^{\max}(M_m) = D_k^{\max}(M_m),$ since the set $D_k^{\max}(M_m)$ is
closed.  Thus, there exists an integer $q,$ a scalar matrix $A = (a_{l,j})
\in M_{q,k},$ and positive matrices, $P_1,\dots,P_q \in M_m^+,$ such
that $\phi(E_{i,j}) = \sum_{l=1}^q \overline{a_{l,i}}P_l a_{l,j}$ for $1
\le i,j \le k.$ If we define positive linear functionals $s_l: M_k
\to \bb C,$ by $s_l(X) = \sum_{i,j} \overline{a_{l,i}} x_{i,j}
a_{l,j},$ then we have that $\phi(E_{i,j}) = \sum_{l=1}^q s_l(E_{i,j})P_l,$
for all $1 \le i,j \le k,$, and hence $\phi(X) = \sum_{l=1}^q s_l(X)P_l$
for every $X \in M_k.$

Conversely, given any positive linear functional $s:M_k \to \bb C$,
we may write $s$ as a sum of vector states; that is, functionals of the form
$X \to \sum_{i,j} \overline{a_i}x_{i,j}a_j.$  Thus if $\phi(X)
= \sum_{l=1}^q s_l(X)P_l$ with each $s_l$ a positive linear functional,
then by increasing the number of terms in the sum we may assume
that each state $s_l$ has the form $s_l(X) = \sum_{i,j}
\overline{a_{l,i}}x_{i,j} a_{l,j},$, and hence $\phi(E_{i,j}) =
\sum_l \overline{a_{l,i}}P_l a_{l,j}$.  Thus $(\phi(E_{i,j})) \in
D_k^{\max}(M_m),$ and it follows that $\phi: M_k \to \OMAX(M_m)$ is
completely positive.
\end{proof}

\begin{cor}\label{c_cpenb}
If $\phi:M_k \to \OMAX(M_m)$ is completely positive, then $\phi$ is entanglement breaking.
\end{cor}
\begin{proof}
By Proposition \ref{pp},
there exist positive linear functionals
$s_l: M_k \to \bb C$ and matrices
$P_l \in M_m^+,$  $l = 1,\dots,q$, such that $\phi(X) = \sum_{l=1}^q s_l(X)P_l.$
Given any $n \in \bb N$ and any positive linear functional
$f:M_n \otimes M_m \rightarrow\bb{C}$, we have that
$$f \circ \phi_n(A \otimes X) =
f(A \otimes \phi(X)) = \sum_{l=1}^q s_l(X) f(A \otimes P_l) = \sum_{l=1}^q s_l(X) g_l(A),$$
where $g_l : M_n \to \bb C$ is the positive linear functional given by $g_l(A) = f(A \otimes P_l)$,
$l = 1,\dots,q$.
Thus $f \circ \phi_n$ is separable.
\end{proof}

\begin{thm}\label{entbr}
Let $\phi:M_k \to M_m$ be a linear map.  Then the following are equivalent:
\begin{itemize}
\item[(i)] $\phi:\OMIN(M_k) \to M_m$ is completely positive.
\item[(ii)] $\phi$ is entanglement breaking.
\item[(iii)] $\phi$ is co-entanglement breaking.
\item[(iv)]  $\phi: M_k \to \OMAX(M_m)$ is completely positive.
\item[(v)] There exist positive linear functionals
$s_l:M_k \to \bb C$ and positive matrices $P_l \in M_m,$ for $1 \leq l \leq K$ such that
$\phi(X) = \sum_{l=1}^K s_l(X)P_l$.
\item[(vi)] There exist rank one matrices $A_l \in M_{k,m}$ for $1 \leq l \leq N$
such that $\phi(X) = \sum_{l=1}^N A_l^*XA_l$.
\item[(vii)] $\phi:\OMIN(M_k) \to \OMAX(M_m)$ is completely positive.
\end{itemize}
\end{thm}
\begin{proof}
The equivalence of (i) and (ii) is stated in Theorem~\ref{th_entb},
while the equivalence of
(iii), (iv), and (v) follows from Corollary~\ref{c_mflat} and Proposition~\ref{pp}.
By Corollary~\ref{c_cpenb}, (iv) implies (ii).
Hence (iii) implies (ii).
Thus, if $\phi = (\phi^{\flat})^{\flat}$ is entanglement breaking,
then $\phi^{\flat}$ is entanglement breaking, and so
(ii) implies (iii).  We now have the equivalence of (i)--(v).

To show that (v) implies (vi), we may assume as before that
each $s_l$ is a vector state. Also, by introducing extra terms we
may assume that each $P_l$ is a rank one positive matrix.  Hence there
exist $V_l \in M_{k,1},$ so that $s_l(X) = V_l^*XV_l$, and there
exist $W_l \in M_{1,m},$ such that $P_l = W_l^*W_l$.  Then (vi) follows by setting $A_l = V_lW_l$.

To see that (vi) implies (v), factor each rank one $A_l = V_lW_l,$
with $V_l \in M_{k,1}$ and $W_l \in M_{1,m}$ and set $s_l(X) =
V_l^*XV_l,$ and $P_l = W_l^*W_l$.

Finally, (vii) clearly implies (i). We show that (v) implies (vii). It
will be enough to consider the case that the sum has a single term,
$\phi(X) = s(X)P.$ But given any operator system $\cl S$, any $P \in \cl
S^+$, and any $s: M_k \to \bb C,$ we see that since $s:\OMIN(M_k) \to \bb C$ is
completely positive, the map $\phi(X) = s(X)P$ is completely
positive.
\end{proof}

We note that the equivalence of (ii) and (vi) of Theorem \ref{entbr}
was proved in \cite{HS}. The following result gives another proof of
the equivalence of (i) and (vii).

\begin{prop} Let $(V,V^+,e)$ be an AOU space and let $\cl S$ be an
  injective operator system. Then $\phi: \OMIN (V) \to \cl S$ is
  completely positive if and only if $\phi: \OMIN (V) \to \OMAX (\cl S)$ is
  completely positive.
\end{prop}
\begin{proof} Assume that $\phi: \OMIN(V) \to \cl S$ is completely
  positive. We have that $\OMIN(V) \subseteq C(S)$ completely order
  isomorphically for some compact, Hausdorff space $S.$ Since $\cl S$
  is injective, we may extend $\phi$ to a completely positive map
  $\psi:C(S) \to \cl S.$ By Proposition~\ref{cxominv}, $C(S)=
  \OMAX(C(S))$ completely order isomorphically, and hence the
  positive map $\psi: C(S) \to \OMAX(\cl S)$ is completely positive.
Hence, $\phi,$ which is the restriction of $\psi$ to $\OMIN(V)$ is
also completely positive.

Conversely, if $\phi: \OMIN(V) \to \OMAX(\cl S)$ is completely
positive, then since the identity map from $\OMAX(\cl S)$ to $\cl S$
is completely positive, we have that $\phi$ is completely positive as
a map into $\cl S.$
\end{proof}

Many of the results about entanglement breaking maps between matrix
algebras can be generalized to operator systems.
The main difference is that the set of separable states need not be
weak*-closed. For this reason, given an operator system $\cl S$ we
call a positive linear functional $\phi:M_n(\cl S) \to \bb C$ {\bf
weak*-separable} if it is a weak*-limit of sums of tensors of
positive linear functionals.

Let $\cl S$ and $\cl T$ be operator systems and let $\phi: \cl S \to \cl
T.$ We call $\phi$ {\bf weak*-entanglement breaking,} if for every
$n,$ and every p.l.f. $s:M_n(\cl T) \to \bb C$ we have that $s \circ
\phi_n: M_n(\cl S) \to \bb C$ is weak*-separable.

The following argument is standard.

\begin{lemma}\label{l_topo}
Let $X$ be a compact Hausdorff space and let
  $\phi:M_n(C(X)) \to \bb C$ be a positive linear functional. Then
  $\phi$ is weak*-separable.
\end{lemma}
\begin{proof}
If $\phi:M_n(C(X)) \to \bb C$ is positive, since $C(X)= \OMIN(C(X))$
completely order isomorphically
by Proposition~\ref{cxominv}, $\phi:
M_n(\OMIN(C(X)))  \to \bb{C}$
is positive. By Remark~\ref{r_wst}, $\phi$ is in the weak*-closure of
$Q_n^{\min}.$ By Remark~\ref{altqnmin}, we see that every element of
$Q_n^{\min}$ is separable.
\end{proof}

\begin{remark} Although every p.l.f. on $M_n(C(X))$ is weak*-separable, in general, a p.l.f. need not be
  separable. An example of a positive linear functional that is not
  separable can be constructed along the same lines as the
  counterexample in Remark~\ref{Dnotarch}.

If we define
$\phi:M_2(C([0,1])) \to \bb C,$ by setting
\[
\phi((f_{i,j})) = \int_0^1 (f_{1,1}(t) + f_{1,2}(t) e^{2 \pi it} +
f_{2,1}(t) e^{-2 \pi it} + f_{2,2}(t)) dt,
\]
then it is fairly easy to see that $\phi$ is a positive linear functional.
A calculation shows that this positive linear functional is not separable.
We sketch the ideas here: First, one shows that if
$\phi = \sum_{j=1}^n s_j \otimes \psi_j,$ where the
$s_j: M_2 \to \bb C$ and $\psi_j: C([0,1]) \to \bb C$
are all positive linear functionals, then the measures associated with
the maps $\psi_j$ are all absolutely continuous with respect to Lebesgue measure,
so that $\psi_j(f) = \int_0^1 f(t) g_j(t) dt$ for some positive Borel functions, $g_j, j=1,\dots,n.$
Next, letting $P_j \in M_2^+$ denote the positive density matrix associated with $s_j, j=1,\dots,n,$ one shows that
$$\begin{pmatrix} 1 & e^{2 \pi it} \\ e^{-2 \pi it} & 1 \end{pmatrix} = \sum_{j=1}^n P_j \otimes g_j(t) \qquad \text{ a.e.}$$
and argues similarly to Remark~\ref{Dnotarch} to show that this is impossible.
\end{remark}

\begin{lemma}\label{l_ominsep}
Let $(V,V^+,e)$ be an AOU space and let $\phi:
  M_n(\OMIN(V)) \to \bb C$ be a positive linear functional.  Then
  $\phi$ is weak*-separable.
\end{lemma}
\begin{proof} We may regard $\OMIN(V)$ as an operator subsystem of
  $C(X)$ where $X$ is the state space of $V.$ We may then extend
  $\phi$ to a positive linear functional on $M_n(C(X))$ and apply
  Lemma \ref{l_topo}.
\end{proof}

\begin{thm} A linear map $\phi: \cl S \to \cl T$ is weak*-entanglement breaking if
  and only if $\phi: \OMIN(\cl S) \to \cl T$ is completely positive.
\end{thm}
\begin{proof}  If $\phi: \OMIN(\cl S) \to \cl T$ is completely positive
  and $s: M_n(\cl T) \to \bb C$ is a positive linear functional, then
  $s \circ \phi_n: M_n(\OMIN(\cl S)) \to \bb C$ is a positive
  linear functional, and hence $s \circ \phi_n$ is weak*-separable by Lemma \ref{l_ominsep}.

Conversely, assume that $\phi: \cl S \to \cl T$ is
weak*-entanglement breaking and let $s:M_n(\cl T) \to \bb C$ be a
positive linear functional.  If $s \circ \phi_n = g \otimes h,$
where $g: \cl S \to \bb C$ and $h: M_n \to \bb C$ are both positive
linear functionals, then let $h = \sum_{i,j} p_{i,j} \delta_{i,j}$
where $P= (p_{i,j}) \in M_n^+.$  Writing $P= X^*X,$ we see that $g
\otimes h = X^* {\rm diag}(g,\dots,g)X \in Q_n^{\min}(\cl S).$ Since
$Q_n^{\min}$ is a cone, every sum of such elementary tensors is in
$Q_n^{\min}.$  Finally, since $s \circ \phi_n$ is weak*-separable
it is in the weak*-closure of $Q_n^{\min}.$ By Remark~\ref{r_wst},
$s \circ \phi_n$ is in $(C_n^{\min})^d.$  Hence the functional $s \circ \phi_n$ is
positive on $C_n^{\min}(\cl S) = M_n(\OMIN(\cl S))^+.$
It follows that $(\phi')_n : M_n(\cl T')\rightarrow M_n({\rm OMIN}(\cl S)')$
is positive for each $n\in \bb{N}$, and hence $\phi : {\rm OMIN}(\cl S)\rightarrow\cl T$
is completely positive.
\end{proof}

Given two normed spaces $X$ and $Y$, the maps in
$CB(\MIN(X),\MAX(Y))$ have been characterized as maps from $X$ to
$Y$ that had factorizations through a Hilbert space that were
``bounded'' and ``co-bounded'' in a certain sense. Consequently,
given two AOU spaces, $V$ and $W$, it is natural to ask for a
characterization of the linear maps that are completely positive
from ${\rm OMIN}(V)$ to ${\rm OMAX}(W).$

\begin{prob}
Given an operator system $\cl S$ and an AOU space $W,$ characterize
the completely
positive maps from $\cl S$ to $\OMAX(W)$. Is such a map a ``limit'' of sums of maps of the
form $\phi(X) = s(X)P$ where $s$ is a positive linear functional on
$\cl S$ and $P \in W^+$?
\end{prob}

\begin{prob}  Given AOU spaces $V$ and $W$, characterize the completely
positive maps from $\OMIN(V)$ to $\OMAX(W).$
\end{prob}

\noindent {\bf Acknowledgements.} The characterizations of
$\OMIN(M_k)$ and $\OMAX(M_k)$ in terms of
entanglement breaking maps given in
the last statements of
Theorem~\ref{th_entb} and Corollary~\ref{c_mflat} were suggested to us
by the referee.

\end{document}